\documentclass[a4paper,11pt]{article}

\usepackage{amsfonts}
\usepackage{amssymb,amsmath,amsthm}
\usepackage{mathtools}

\usepackage{authblk}
\usepackage[numbers]{natbib}

\usepackage[Algorithm]{algorithm}
\usepackage{algpseudocode}
\usepackage{placeins}
\let\Emptyset\emptyset
\usepackage{graphicx}        
\usepackage{multicol}        
\usepackage{mdframed}
\usepackage{tikz}
\usetikzlibrary{shapes,arrows}
\usetikzlibrary{patterns}

\usepackage{pgfplots}
\pgfplotsset{compat=1.13}
\usepgfplotslibrary{groupplots,fillbetween}

\definecolor{HSUred}{RGB}{197,0,66}

\usepackage{siunitx}
\usepackage{hyperref}
\usepackage{xcolor}
\usepackage{subcaption}
\usepackage{graphicx}
\usepackage{booktabs}
\usepackage{float}

\renewcommand{\emptyset}{\Emptyset}

\newcommand{\velocity}{\boldsymbol v}
\newcommand{\pressure}{p}
\newcommand{\normal}{\boldsymbol n}
\newcommand{\primalu}{\boldsymbol{u}}
\newcommand{\primalphi}{\boldsymbol{\varphi}}
\newcommand{\primalpsi}{\boldsymbol{\psi}}
\newcommand{\dualz}{\boldsymbol{z}}
\newcommand{\dualvelocity}{\boldsymbol w}
\newcommand{\dualpressure}{q}
\newcommand\tsb[1]{_{\textup{#1}}}
\newcommand\Ieff{I\tsb{eff}}
\newcommand{\mc}[1]{\multicolumn{1}{c}{#1}}

\newtheorem{defi}{Definition}[section]
\newtheorem{theorem}[defi]{Theorem}

\newtheorem{remark}[defi]{Remark}

\newenvironment{mproof}{\paragraph{Proof.}}{\hfill$\blacksquare$}
\numberwithin{equation}{section}
\numberwithin{table}{section}
\numberwithin{figure}{section}

\usepackage{geometry}
\begin{document}

\title{\Large \textbf{Goal-Oriented Space-Time Adaptivity for the Navier--Stokes Equations based on the Dual Weighted Residual Method}}
\author[M.\ P.\ Bruchh\"auser, N.\ Margenberg, M.\ Bause]
{\large \textbf{Marius Paul Bruchh\"auser}\thanks{$^{,\ddagger}$ Helmut Schmidt University Hamburg, Faculty of Mechanical and Civil Engineering, Chair of Numerical Mathematics,
Holstenhofweg 85, 22043 Hamburg, Germany,\\ bruchhaeuser@hsu-hamburg.de ($^\ast$corresponding author),
bause@hsu-hamburg.de} $^{,\ddagger}$
$\,\boldsymbol{\cdot}$
\textbf{Nils Margenberg}\thanks{University of Magdeburg, Institute for Analysis and Numerics, Universit\"atsplatz 2, 39104 Magdeburg, Germany, nils.margenberg@ovgu.de.de}
$\,\boldsymbol{\cdot}$
\textbf{Markus Bause$^{\ddagger}$}\\
}
\date{}
\maketitle

\begin{abstract}
\noindent
This work presents a goal-oriented a posteriori error estimator based on the Dual Weighted Residual (DWR) method together with space-time mesh adaptivity for the Navier--Stokes equations. The resulting nonlinear algebraic systems on the space-time slabs are solved by Newton's method with GMRES, preconditioned by a slab-wise geometric multigrid method. This combination yields reliable control of target quantities on computationally feasible space-time meshes together with a robust and efficient solution of the algebraic systems. The implementation is based on a MPI-parallel programming model in the deal.II library. Further ingredients are a discontinuous Galerkin discretization in time and inf-sup stable finite element pairs with discontinuous pressure on tensor-product meshes. The performance of the approach is investigated in benchmark computations with regard to accuracy, efficiency, and stability.
\end{abstract}


\section{Introduction}
\label{sec:1:introduction}
The accurate numerical simulation of incompressible viscous fluid flow governed by the Navier--Stokes equations remains challenging, in particular when strong dynamics and fine-scale structures in space and time have to be resolved. The robust and efficient solution of the arising algebraic systems adds a further layer of complexity. Moreover, algorithms that address these challenges have to be implemented in a software framework that is suitable for computations on modern hardware architectures with tens of thousands of compute cores, including hybrid CPU and GPU systems. This places high demands on the underlying data structures. In this work, we design an approach that addresses these requirements. The implementation is based on the deal.II library \cite{BruchhaeuserABBFGHHKKMMPPSTWZ24}.

Space-time mesh adaptation based on rigorous error representations, or on indicators derived from such representations, has been shown to yield accurate results on computationally feasible grids at moderate numerical cost, even for problems with complex solution structures \cite{langer-2021,SteinbachYang2019}. Among these approaches, the Dual Weighted Residual method (DWR) \cite{BruchhaeuserBR01} offers the possibility of controlling the discretization error in user-defined goal quantities of practical interest. The key idea of this approach is to compute dual weights for the local residuals entering the error representation. This requires the solution of the generally nonlinear primal problem and its linearized dual problem by space-time finite element methods. The DWR method has been successfully applied to numerous problem classes, including fluid flow \cite{BruchhaeuserBR12,RePEc:spr:978-3-642-18775}, wave propagation \cite{BruchhaeuserBGR10}, fluid-structure interaction \cite{BruchhaeuserR17,BruchhaeuserWW21}, among others.
For a comprehensive review with detailed references on various model problems, we refer, for instance, to  \cite{BruchhaeuserBR01,BruchhaeuserBR03,BruchhaeuserR17}.
Error control in norms associated with the weak formulation constitutes a conceptually different class of techniques for a posteriori error estimation and automatic mesh adaptation. For this, an approach using equilibrated flux reconstructions is proposed in \cite{ern-2017}.

The DWR concept relies on solving the primal Navier--Stokes system and the associated linearized dual problem by space-time finite element methods. Both, lower- and higher-order approximations, in space and time are of interest. They should maintain stability and preserve as much of the structure of the continuous problem as possible. For the time discretization we use the discontinuous Galerkin method \cite{10.5555/1211262,hussain-2012}. In space, we employ inf-sup stable finite element pairs with discontinuous pressure. Discontinuous pressure approximations have demonstrated advantages over continuous ones  \cite{10.1002/fld.195}. Moreover, they are advantageous for the construction and application of the slab-wise local Vanka smoother within the geometric multigrid preconditioner. Nonhomogeneous Dirichlet boundary conditions are imposed by Nitsche's method~\cite{BruchhaeuserN71}, which leads to a strictly variational formulation rather than enforcing the boundary data through the definition of the function spaces and algebraic systems.

We localize the space-time finite element approximation of the primal and dual DWR problems to space-time slabs $\Omega\times (t_{n-1},t_n]$, where $\Omega \subset \mathbb{R}^d$ with $d \in \{2,3\}$. This yields a time-marching approach. For higher polynomial degrees in time and high spatial resolution, the arising algebraic systems already become large, and the Newton-linearized Jacobian has a complex block structure that complicates the construction of efficient preconditioners. Holistic space-time finite element frameworks, in which all time steps are solved fully coupled and simultaneously, have also been studied yet; cf., e.g., \cite{10.1137/18M1172466,10.1515/9783110548488,WIENERS2023294,falgout-2017}. Their feasibility continues to remain an open problem if three space variables and a huge number of time steps are involved. The lack of efficient preconditioning techniques for such systems still limits the application of holistic approaches. For this reason, we use a slab-wise approach with a tensor-product mesh structure in this work.

For the slab-wise solution of the Newton-linearized primal problems and of the linear dual problems, we apply a Flexible Generalized Minimum RESidual (FGMRES) method preconditioned by a single $V$-cycle of a geometric multigrid (GMG) method. Nowadays, multigrid techniques are typically used as preconditioners for Krylov subspace methods rather than as stand-alone solvers. GMG methods have been studied extensively for the iterative solution of discrete systems arising from stationary partial differential equations. For this we refer to the pioneering work \cite{hackbusch1985multigrid}. Owing to their optimal linear complexity, GMG methods are highly efficient and robust for this class of problems \cite{gmeiner-2015}. They have also been generalized in several directions, including $hp$ multigrid methods based on hierarchies of polynomial degrees and meshes. For space-time finite element approximations of the Navier--Stokes equations, block systems with generalized saddle-point structure arise, where the precise block structure depends on the ordering of the slab unknowns \cite{BruchhaeuserAB23,25M1734142}. For such systems, the smoothing procedure is decisive for the overall efficiency of the multigrid preconditioner. In the present work, we use a local Vanka smoother, which can be interpreted as a non-overlapping Schwarz method. It is among the most effective smoothers for saddle-point problems, although its cost increases with the polynomial degree and with the number of space-time unknowns \cite{25M1734142}. Algebraic MultiGrid (AMG) methods, a further class of multigrid techniques, perform damping of high frequency errors in the solution to linear systems and coarse grid corrections on the algebraic level and do not need any geometric information for building a hierarchy of problems \cite{10.1137/1.9781611971057,steinbach-2018}. However, memory usage limits the application of AMG if large scale systems are encountered.

In this work, we study the combination of DWR error estimation, automatic space-time mesh adaptation, and slab-wise geometric multigrid preconditioning for the Navier--Stokes equations. The focus is on goal-oriented adaptivity and its interaction with a GMG-preconditioned Newton--Krylov method. For this purpose, we combine a DWR a posteriori error estimator with adaptive space-time refinement and GMRES iterations preconditioned by slab-wise geometric multigrid. The implementation uses MPI parallelization in deal.II, discontinuous Galerkin time discretization, and inf-sup stable finite element pairs with discontinuous pressure on tensor-product meshes. On benchmark computations, we assess estimator quality, reduction of the goal error relative to total computational work, and solver performance for both primal and adjoint problems. A central question we address is whether the gain from goal-oriented adaptivity justifies the additional cost of the adjoint solve and how this trade-off interacts with the robustness of the solvers.

This work is organized as follows.
In Sec.~\ref{sec:2:weak}, we present the weak form as well as the space-time discretization schemes.
In Sec.~\ref{sec:3:error}, we derive an a posteriori error representation for the model problem.
In Sec.~\ref{sec:4:practical}, we discuss practical aspects of the implementation as well as the underlying space-time adaptive algorithm.
In Sec.~\ref{sec:5:numerical}, we present the results of our numerical examples.
Finally, Sec.~\ref{sec:6:conclusion} summarizes the results and outlines directions for future research.

\section{Model Problem and Discretization in Space and Time}
\label{sec:2:weak}
In this work, we investigate the accurate and efficient numerical simulation of incompressible viscous flow modeled by the nonstationary Navier--Stokes equations
\begin{equation}
\label{eq:1:navier}
\begin{array}{rcl @{\,\,}l @{\,\,}l @{\,}l}
\partial_t \velocity
- \nu \Delta \velocity
+\velocity\cdot\nabla\velocity
+ \nabla \pressure
&=& \boldsymbol{f}
& \textnormal{in} & Q & = \Omega \times I\,,\\[.5ex]
\nabla \cdot \velocity &=& 0
& \textnormal{in} & Q & = \Omega \times I\,,\\[.5ex]
%
%
\velocity &=& \velocity_D
& \textnormal{on} & \Sigma_{D}
& = \Gamma_{D} \times I\,,\\[.5ex]
%
%
\nu\partial_{\normal}\velocity - \pressure\,\normal
&=& \boldsymbol 0 & \textnormal{on} & \Sigma_{\textnormal{outflow}}
& = \Gamma_{\textnormal{outflow}} \times I\,,\\[.5ex]
\velocity(0) &=& \velocity_0
& \textnormal{on} & \Sigma_0 & = \Omega \times \{ 0 \}\,,\\[.5ex]
\end{array}
\end{equation}
In \eqref{eq:1:navier}, we denote by $Q=\Omega \times I$ the space-time domain, where $\Omega\subset \mathbb{R}^d$, with $d=2$ or $d=3$, is a polygonal or polyhedral bounded domain with Lipschitz boundary $\partial\Omega$ and $I=(0,T)$, with $0 < T < \infty$,
is a finite time interval.
The flow behavior at the boundary $\partial\Omega = \Gamma_{D}
\cup \Gamma_{\textnormal{outflow}}$ is modeled by a function $\velocity_D$
on the segments $\Gamma_D\coloneq\Gamma_{\textnormal{inflow}} \cup \Gamma_{\textnormal{wall}}$,
that prescribes an inflow profile $\velocity_{\textnormal{in}}$ on
$\Gamma_{\textnormal{inflow}}$ and a no-slip condition $\velocity=\boldsymbol 0$ on the fixed wall segment $\Gamma_{\textnormal{wall}}$, as well as a so-called ``do-nothing'' condition on the artificial boundary segment $\Gamma_{\textnormal{outflow}}$. For this refer to \cite{BruchhaeuserJ16,BruchhaeuserHRT96}.
The velocity field $\velocity$ and the pressure variable $\pressure$ are the unknown variables. The parameter $\nu>0$ denotes the fluid's viscosity,
$\boldsymbol{f}$ is a given volume force density, and
$\velocity_{0}$ is the initial value.

For the sake of brevity, we present the continuous and semidiscrete in time variational problem for homogeneous Dirichlet boundary conditions only. The mixed boundary conditions in \eqref{eq:1:navier} are incorporated in the fully discrete system and derivation of the duality based a posteriori error estimation; cf.~Remark~\ref{Remark:Forms}. Further, they are used in one of our numerical experiments. Throughout, we use standard notation for Sobolev and Bochner spaces. In particular, by $L_\sigma^2(\Omega)^d$ we denote the $L^2$-closure of the set of divergence-free test functions on $\Omega$. For the existence, uniqueness and regularity of weak solutions to the Navier--Stokes equations such that all of the arguments and terms used below are well-defined, we refer to the literature, cf., e.g., \cite{BruchhaeuserHR82,BruchhaeuserJ16,1977iii}.
\subsection{Weak Form and Discretization in Space and Time}
\label{sec:2.1:weak}
For the discretization in space and time we follow the lines of
\cite{BruchhaeuserBR12,BruchhaeuserAB23}.
The variational formulation of \eqref{eq:1:navier} reads as:
\newline
\textit{For $\boldsymbol{f} \in L^2(I;H^{-1}(\Omega)^{d})$ and
$\velocity_{0}\in L_\sigma^2(\Omega)^d,$\; find $\primalu\coloneq\{\velocity,\pressure\}\in X$
such that}
\begin{equation}
\label{eq:2:weakform}
A(\primalu)(\primalphi) = F(\primalpsi) \quad \forall \primalphi\coloneq\{\primalpsi,\xi\} \in X\,,
\end{equation}
\textit{where the semi-linear form
$A(\cdot)(\cdot)$
and the linear form
$F(\cdot)$
are defined by}
\begin{eqnarray}
\label{eq:3:SLFA}
A(\primalu)(\primalphi) & \coloneq &
\int_{I}\big\{(\partial_t \velocity,\primalpsi)
+ a(\primalu)(\primalphi)
\big\} \mathrm{d} t\;
+ (\velocity(0),\primalpsi(0) )\,,
\\[1.5ex]
\label{eq:4:LiFormF}
F(\primalpsi) & \coloneq & \int_I(\boldsymbol{f},\primalpsi)\;\mathrm{d}t + (\velocity_{0},\primalpsi(0) ) \,,
\end{eqnarray}
\textit{with the inner semi-linear form
$a(\cdot)(\cdot)$
given by}
\begin{equation}
\label{eq:5:inSLFa}
a(\primalu)(\primalphi) \coloneq \nu(\nabla \velocity, \nabla \primalpsi)
+(\velocity\cdot \nabla \velocity, \primalpsi)
-(\pressure,\nabla\cdot\primalpsi)
+ (\nabla\cdot\velocity,\xi)\,.
\end{equation}
Here, $(\cdot,\cdot)$ denotes the inner product on $L^2(\Omega)$ (or $L^2(\Omega)^d$)
and the underlying function space $X$ is defined as
\begin{equation}
\label{eq:6:X}
X\coloneq\{\primalu=\{\velocity,\pressure\}, \velocity\in L^2\big(I;H_0^1(\Omega)^d\big),
\partial_t\velocity\in L^2\big(I;H^{-1}(\Omega)^d\big), \pressure\in L^2\big(I;L_0^2(\Omega)\big)\}\,.
\end{equation}
For the semi-discretization in time we use a discontinuous Galerkin method dG($r$)
with an arbitrary polynomial degree $r\ge 0$.
Let $\mathcal T_{\tau}\,=\{0\}\cup I_1 \cup \cdots \cup I_N$ be a partition of the closure of the time domain
$\bar{I}=[0,T]$ with left-open subintervals $I_n\coloneq(t_{n-1},t_n]$, $n=1,\dots,N$, with $0=:t_0<t_1<\dots<t_N\coloneq T$ and time step sizes $\tau_n=t_n-t_{n-1}$ and the global time discretization parameter $\tau=\max_{n}\,\tau_{n}$.
On the subintervals $I_n$, we define the time-discrete function space $X_{\tau}^{r}$ as
\begin{equation}
\label{eq:7:Def_X_tau_r}
\begin{aligned}
X_{\tau}^{r} \coloneq
 \Big\{
 & \primalu_{\tau}=\{\velocity_{\tau},\pressure_{\tau}\}\in L^{2}\big(I; H_0^1(\Omega)^d\times L_0^2(\Omega)\big)\big|
 \velocity_{\tau}|_{I_{n}}\in \mathcal{P}_{r}\big(I_{n}; H_0^1(\Omega)^d\big),\\
 & \velocity_{\tau}(0)\in L^2(\Omega)^d, \pressure_{\tau} \in \mathcal{P}_{r}\big(I_{n}; L_0^2(\Omega)\big)\,, n=1,\dots,N
\Big\}\,,
\end{aligned}
\end{equation}
where $\mathcal{P}_{r}(I_{n}; Y)$ denotes the space of all
polynomials in time up to degree $r\geq0$ on $I_n$ with values in a Banach space $Y$.
For some discontinuous in time function
$\velocity_{\tau}\in X_{\tau}^{r}$ we define
the limits $\velocity_{\tau,n}^{\pm}$ from above and below of
$\velocity_{\tau}$ at $t_n$ as
well as their jump at $t_n$ by
\begin{displaymath}
\begin{array}{lcrclcr}
\velocity_{\tau,n}^{\pm}
& \coloneq &
\displaystyle\lim_{t\mapsto t_n\pm0} \velocity_\tau(t) \,,
&
[\velocity_{\tau}]_{n} & \coloneq & \velocity_{\tau,n}^{+}
-\velocity_{\tau,n}^{-} \,.
\end{array}
\end{displaymath}
Then, the semi-discretization in time scheme of the incompressible Navier--Stokes equations \eqref{eq:2:weakform} reads as follows:
\textit{Find $\primalu_{\tau}=\{\velocity_{\tau},\pressure_{\tau}\} \in X_{\tau}^{r}$
such that
}
\begin{equation}
\label{eq:8:time-discrete}
A_{\tau}(\primalu_{\tau})(\primalphi_{\tau})
=
F_\tau(\primalpsi_{\tau})
\quad \forall \primalphi_{\tau}=\{\primalpsi_{\tau},\xi_{\tau}\} \in X_{\tau}^{r}\,,
\end{equation}
\textit{where the semi-discrete semi-linear form $A_{\tau}(\cdot)(\cdot)$ and linear form $F_{\tau}(\cdot)$ are defined by}
\begin{subequations}
\label{eq:9:Def_A_tau_Def_F_tau}
\begin{align}
\label{eq:9:Def_A_tau}
A_{\tau}(\primalu_{\tau})(\primalphi_{\tau})
& \coloneq
\sum_{n=1}^{N}\int_{I_n}\big\{(\partial_{t} \velocity_{\tau},\primalpsi_{\tau})
+ a(\primalu_{\tau})(\primalphi_{\tau})
\big\} \mathrm{d} t\;
\\
\nonumber
& + \sum_{n=2}^{N}([\velocity_\tau]_{n-1},\primalpsi_{\tau,n-1}^+ )
+ (\velocity_{\tau,0}^+,\primalpsi_{\tau,0}^+)\,,
\\
\label{eq:9:Def_F_tau}
F_\tau (\primalpsi_{\tau}) & \coloneq  \int_I(\boldsymbol{f},\primalpsi_{\tau})\;\mathrm{d}t
+ (\velocity_{0},\primalpsi_{\tau,0}^+)\,.
\end{align}
\end{subequations}

\begin{remark}
By the choice of a discontinuous test basis, supported on the subintervals $I_n$, Equation \eqref{eq:8:time-discrete} splits into a sequence of variational problems on the space-time slabs $\Omega\times I_n$. This allows us to solve the algebraic problems associated with the slabs sequentially in a time-marching approach. Further, we introduce dynamic spatial meshes below, permitting the spatial decomposition of $\Omega$ to depend on the subintervals $I_n$ and change from subinterval to subinterval. For the analysis and application of dynamic meshes we refer to \cite{10.1002/fld.2625,10.1093/imanum/drt001,BruchhaeuserSV08}.
\end{remark}

For the discretization in space, let $\mathcal{T}_h^{n}=\{K\}$ be a triangulation of mesh size $h_n$ of the domain $\Omega$ on subinterval $I_n$ such that $\overline{\Omega}=\bigcup_{K\in\mathcal{T}_h^{n}}\overline{K}$. For $t=0$, we let $\mathcal T_h^0$ denote a decomposition of $\Omega$ with mesh size $h_0$. The maximal mesh size is denoted by $h\coloneq \max_{n} h_n$. We define the space-time mesh slabs $Q_n \coloneq \mathcal{T}_h^n \times I_n$ for $n=1,\dots, N$. Our mesh adaptation process yields locally refined elements, which is enabled by using hanging nodes. We point out that the global conformity of the finite element approach is preserved since the unknowns at such hanging nodes are eliminated by interpolation between the neighboring 'regular' nodes; cf.~\cite[Chapter 4.2]{BruchhaeuserBR03} and \cite{BruchhaeuserCO84}.
%

%
In our GMG approach we apply a local (element-based) Vanka smoother \cite{BruchhaeuserAB23}. Numerical studies have shown, that this smoother performs efficiently for local velocity-pressure couplings by using discontinuous pressure approximations. Therefore, we apply the inf-sup stable family of finite element pairs $Q_k/P_{k-1}^{\textnormal{disc}}$, for $k\geq 2$ \cite{BruchhaeuserJ16,BruchhaeuserMT02}). Thus, we built our spatial discretization on the following pair of spaces:
\begin{equation}
  \label{eq:10:VhQh}
  \begin{aligned}
  \boldsymbol V_h^{k}(\mathcal{T}_h^n) &\coloneq \{\boldsymbol v_h\in [H^1(\Omega)]^{d}: (\boldsymbol v_h)_{|K}\in
  V^{k}(K)\
  \forall K\in\mathcal T_h^n\}\,,\\
  Q_h^{k-1}(\mathcal{T}_h^n) &\coloneq \{ q_h\in L_0^2(\Omega): (q_h)_{|K}\in
  Q^{k-1}(K)\
  \forall K\in\mathcal T_h^n\}\,,
  \end{aligned}
\end{equation}
where the finite elements $\boldsymbol V^{k}(K)$ and $Q^{k-1}(K)$ are defined by
\begin{equation}
\label{eq:11:VH}
\boldsymbol V^{k}(K)\coloneq [Q_k]^d\circ F_K^{-1}\,,\qquad
 Q^{k-1}(K)\coloneq P_{k-1}^{\textnormal{disc}}\circ F_K^{-1}\,,
\end{equation}
and $F_K$ is the standard multilinear map from the reference element to $K$. We employ the mapped variant of $P_{k-1}^{\mathrm{disc}}$ for geometric consistency on curved and non-affine meshes and improved conditioning of the algebraic systems.
Here, ${Q}_k, k\geq 2$, denotes the space of polynomials defined on the reference element
with maximum degree $k$ in each variable and ${P}_{k-1}^{\textnormal{disc}}$ is the space of polynomials
defined on the reference element with maximum degree $k-1$.
Now, by replacing $H_0^1(\Omega)$ and $L_0^2(\Omega)$ in the
definition of the semi-discrete function space $X_{\tau}^{r}$
in (\ref{eq:7:Def_X_tau_r}) by the spatial finite element spaces defined above,
we obtain the fully discrete function space
\begin{equation}
\begin{aligned}
\label{eq:12:Def_X_tau_h_r_k}
X_{\tau h}^{r,k} \coloneq \Big\{ &
\primalu_{\tau h}=\{\velocity_{\tau h},\pressure_{\tau h}\}\, \big|\,
\velocity_{\tau h}|_{I_n} \in \mathcal{P}_r\big(I_n;
\bold V_{h}^{k}(\mathcal{T}_h^n)
\big)
\,, \velocity_{\tau h}(0) \in
\bold V_{h}^{k}(\mathcal{T}_h^0)
\,,\\
& \pressure_{\tau h}|_{I_n} \in \mathcal{P}_r\big(I_n;
Q_{h}^{k-1}(\mathcal{T}_h^n)
\big)\,,
n=1,\dots,N
\Big\}\,.
\end{aligned}
\end{equation}
For the treatment of Dirichlet boundary conditions, we use Nitsche's method (cf.~\cite{BruchhaeuserN71}), where the Dirichlet boundary conditions are imposed in weak form. For this, we let $N_{\Gamma_D}: L^2\big(I;H^{1/2}(\Gamma_D)^d\big)\times X_{\tau h}^{r,k} \to \mathbb R$ be defined by
\begin{equation*}
N_{\Gamma_D}(\velocity, \primalphi) x\coloneq
\sum_{n=1}^{N}\int_{I_n}\big\{
-(\velocity, \nu\partial_{\normal}\primalpsi+\xi\,\normal)_{\Gamma_D}
+\nu\dfrac{\gamma_1}{h}(\velocity, \primalpsi)_{\Gamma_D}
+\dfrac{\gamma_2}{h}(\velocity\cdot\normal, \primalpsi\cdot\normal)_{\Gamma_D}
\big\}\;\mathrm{d} t
\end{equation*}
for $\velocity \in L^2\big(I;H^{1/2}(\Gamma_D)^d$ and $\primalphi \in X_{\tau h}^{r,k}$ and some numerical (tuning) parameters $\gamma_1>0$ and $\gamma_2>0$ ensuring stability of the discrete system .
Further, we define ${B}_{\Gamma_D}: X_{\tau h}^{r,k} \times X_{\tau h}^{r,k}\to \mathbb R$ by
\begin{equation}
\label{eq:13:nitsche}
{B}_{\Gamma_D}({\primalu},\primalphi)\coloneq
\sum_{n=1}^{N}\int_{I_n}
-(\nu\partial_{\normal}{\velocity}-{\pressure}\,\normal,\primalpsi)_{\Gamma_D}\,\mathrm{d} t
+N_{\Gamma_D}({\velocity}, \primalphi)
\end{equation}
for $(\primalu,\primalphi)\in X_{\tau h}^{r,k} \times X_{\tau h}^{r,k}$.
Finally, the fully discrete problem to Equation~\eqref{eq:2:weakform}
reads as follows: \textit{Find $\primalu_{\tau h}=\{\velocity_{\tau h},\pressure_{\tau h}\} \in X_{\tau h}^{r,k}$
such that
}
\begin{equation}
\label{eq:14:fully-discrete}
A_{\tau h}(\primalu_{\tau h})(\primalphi_{\tau h})
=
F_{\tau h}(\primalphi_{\tau h})
\quad \forall \primalphi_{\tau h}=\{\primalpsi_{\tau h},\xi_{\tau h}\} \in X_{\tau h}^{r,k}\,,
\end{equation}
\textit{where the semi-discrete semi-linear form $A_{\tau h}(\cdot)(\cdot)$ and the linear form $F_{\tau h}(\cdot)$ are defined by}
\begin{subequations}
\label{eq:15:Def_A_tau_h_Def_F_tau_h}
\begin{align}
\label{eq:15:Def_A_tau_h}
A_{\tau h}(\primalu_{\tau h})(\primalphi_{\tau h})
\coloneq &
A_{\tau}(\primalu_{\tau h})(\primalphi_{\tau h})
+{B}_{\Gamma_D}(\primalu_{\tau h},\primalphi_{\tau h})
\\
\label{eq:15:Def_F_tau_h}
F_{\tau h}(\primalphi_{\tau h}) \coloneq &
F_{\tau} (\primalpsi_{\tau h})
+N_{\Gamma_D}(\velocity_D, \primalphi_{\tau h})\,.
\end{align}
\end{subequations}
We note that in the case of homogeneous Dirichlet boundary conditions the last term in~\eqref{eq:15:Def_F_tau_h} vanishes since $\velocity_D=\boldsymbol 0$.

\begin{remark}
\label{Remark:Forms}
We note that in the continuous problem \eqref{eq:2:weakform} and discrete in time system \eqref{eq:8:time-discrete} the homogeneous Dirichlet boundary conditions are prescribed in strong form by the definition of the solution spaces. In the fully discrete problem \eqref{eq:14:fully-discrete}, the Dirichlet boundary conditions are implemented in weak or variational form. This affects our duality-based a posteriori error control in Section~\ref{sec:3:error} below and the derivation of the error identity therein.
This approach follows a more general version compared to the result derived in \cite{BruchhaeuserBR01}, covering the case of inconsistently stabilized discretizations and can be found in \cite[Lem.~5.1]{BruchhaeuserBR12}.
Here, the additional terms are introduced only at the fully discrete level, without corresponding consistent counterparts in the continuous or semi-discrete in time formulations, cf.~\cite{BruchhaeuserBR12,BruchhaeuserB22}.
\end{remark}


\section{DWR-Based A Posteriori Error Estimation}
\label{sec:3:error}
In this section, we derive an a posteriori error representation based on the
DWR method.
This representation is given in terms of a user-chosen goal functional $J:X\mapsto\mathbb{R}$,
in general given as
\begin{equation}
\label{eq:3:1:Jgeneral}
J(\primalu)=\int_0^T J_1(\primalu(t))\;\mathrm{d}t + J_2(\primalu(T))\,,
\end{equation}
where $J_1$ and $J_2$
are three times differentiable functionals and either of them may vanish;
cf., e.g., \cite{BruchhaeuserR17}.
Since we aim at controlling the error due to the discretization in time and in
space, we split the a posteriori error representation with respect to $J$ into
the contributions
\begin{equation}
\label{eq:3:2:Jsplitted}
 J(\primalu)-J(\primalu_{\tau h}) =
 J(\primalu)-J(\primalu_{\tau})
 + J(\primalu_{\tau})-J(\primalu_{\tau h})\,.
\end{equation}
To derive these representations, we introduce the Lagrangian functionals
$\mathcal{L}: X\times X \rightarrow \mathbb{R}$,
$\mathcal{L}_\tau: X_{\tau}^{r} \times X_{\tau}^{r}
\rightarrow \mathbb{R}$, and
$\mathcal{L}_{\tau h}:
X_{\tau h}^{r,k} \times X_{\tau h}^{r,k}
\rightarrow \mathbb{R}$
based on the Euler Lagrange method of constraint optimization:
\begin{subequations}
\label{eq:3:3:Lagrangians}
\begin{align}
\label{eq:3:3:Def_L_u_z}
\mathcal{L}(\primalu,\dualz) & \coloneq  J(\primalu)
+ F(\dualvelocity)
- A(\primalu)(\dualz)\,,
\\
\label{eq:3:3:Def_L_tau_u_z}
\mathcal{L}_{\tau}(\primalu_{\tau},\dualz_{\tau}) & \coloneq
J(\primalu_{\tau}) + F_\tau(\dualvelocity_{\tau})
- A_{\tau}(\primalu_{\tau})(\dualz_{\tau})\,,
\\
\label{eq:3:3:Def_L_tau_h_u_z}
\mathcal{L}_{\tau h}(\primalu_{\tau h},\dualz_{\tau h}) & \coloneq
J(\primalu_{\tau h})
+ F_{\tau h} (\dualz_{\tau h})
- A_{\tau h}(\primalu_{\tau h})(\dualz_{\tau h})
\\
\nonumber
&=\mathcal{L}_{\tau}(\primalu_{\tau h},\dualz_{\tau h})
-{B}_{\Gamma_D}(\primalu_{\tau h},\dualz_{\tau h})
+N_{\Gamma_D}(\velocity_D,\dualz_{\tau h})
\,.
\end{align}
\end{subequations}
Here, the Lagrange multipliers $\dualz=\{\dualvelocity,\dualpressure\}$,
$\dualz_{\tau}=\{\dualvelocity_{\tau},\dualpressure_{\tau}\},$ and
$\dualz_{\tau h}=\{\dualvelocity_{\tau h},\dualpressure_{\tau h}\}$
are called dual variables in contrast to the primal variables
$\primalu=\{\velocity,\pressure\}$, $\primalu_{\tau}=\{\velocity_{\tau},\pressure_{\tau}\},$ and
$\primalu_{\tau h}=\{\velocity_{\tau h},\pressure_{\tau h}\}$; cf. \cite{BruchhaeuserBR01,BruchhaeuserBR12}.

Considering the directional derivatives of the Lagrangian functionals, also
known as G\^{a}teaux derivatives, with respect to their first argument, i.e.\
\begin{equation}
\label{eq:3:4:Def_Gateaux_derivative}
\mathcal{L}^{\prime}_{\primalu}(\primalu,\dualz)(\primalphi) \coloneq
\lim_{t\neq0,t\rightarrow 0}
t^{-1}\big\{\mathcal{L}(\primalu + t\primalphi,\dualz)
-\mathcal{L}(\primalu,\dualz)\big\},
\quad \primalphi \in X\,,
\end{equation}
leads to the so-called dual problems, cf., e.g., \cite{BruchhaeuserBR12}.
The continuous, semi-discrete, and fully discrete dual solutions
$\dualz~\in~X$, $\dualz_\tau~\in~X_{\tau}^{r},$ and
$\dualz_{\tau h} \in X_{\tau h}^{r,k}$
are determined by the optimality conditions
\begin{subequations}
\label{eq:3:5:Def_L_uprime_u_z_Def_L_tau_uprime_u_z_Def_L_tauh_uprime_u_z}
\begin{align}
\label{eq:3:5:Def_L_uprime_u_z}
\mathcal{L}^{\prime}_{\primalu}(\primalu,\dualz)(\primalphi) & =  0
\quad \forall \primalphi \in X\,,
\\
\label{eq:3:5:Def_L_tau_uprime_u_z}
\mathcal{L}^{\prime}_{\tau, \primalu}(\primalu_\tau,\dualz_\tau)(\primalphi_{\tau}) & =  0
\quad \forall \primalphi_{\tau} \in
X_{\tau}^{r}\,,
\\
\label{eq:3:5:Def_L_tauh_uprime_u_z}
\mathcal{L}^{\prime}_{\tau h, \primalu}(\primalu_{\tau h},\dualz_{\tau h})(\primalphi_{\tau h}) & =  0
\quad \forall
\primalphi_{\tau h} \in X_{\tau h}^{r,k}\,.
\end{align}
\end{subequations}
More precisely, the continuous dual solution $\dualz \in X$ is the solution of
\begin{equation}
\label{eq:3:6:dual_continuous}
A^{\prime}(\primalu)(\primalphi,\dualz)
=
J^{\prime}(\primalu)(\primalphi)
\quad \forall \primalphi \in X\,,
\end{equation}
where
the adjoint semi-linear form $A^{\prime}(\cdot)(\cdot, \cdot)$ is given by
\begin{equation}
\label{eq:3:7:Def_A_prime_u_phi_z}
A^{\prime}(\primalu)(\primalphi,\dualz)
\coloneq
\int_I\big\{(\primalpsi,-\partial_{t} \dualvelocity)
+a^{\prime}(\primalu)(\primalphi,\dualz)\big\}\mathrm{d}t
+(\primalpsi(T),\dualvelocity(T))\,.
\end{equation}
We note that for the representation \eqref{eq:3:7:Def_A_prime_u_phi_z} of
$A^{\prime}(\cdot)(\cdot, \cdot)$ integration by parts in time is applied,
which is allowed for weak
solutions $\dualvelocity\in X$; cf., e.g., \cite[Lemma 8.9]{BruchhaeuserR17}.
The derivative $a^{\prime}(\primalu)(\primalphi,\dualz)$ of the inner
semi-linear form $a(\primalu)(\dualz)$ in $A^{\prime}$ admits
the explicit form
\begin{displaymath}
a^{\prime}(\primalu)(\primalphi,\dualz)=
\nu(\nabla \primalpsi,\nabla \dualvelocity)
+ (\velocity \cdot \nabla \primalpsi,\dualvelocity)
+ (\primalpsi\cdot\nabla\velocity,\dualvelocity)
-(\xi,\nabla\cdot\dualvelocity)
+(\nabla\cdot\primalpsi,\dualpressure)\,.
\end{displaymath}
The right-hand side of Eq.~\eqref{eq:3:6:dual_continuous} is
given by
\begin{equation}
\label{eq:3:8:Def_J_prime_u_phi}
J^{\prime}(\primalu)(\primalphi) \coloneq
\int_I J_1^{\prime}(\primalu)(\primalphi)\mathrm{d}t+J_2^{\prime}(\primalu(T))(\primalphi(T))\,.
\end{equation}
Further, the semi-discrete dual solution $\dualz_\tau\in X_{\tau}^{r}$
and the fully discrete dual solution $\dualz_{\tau h} \in X_{\tau h}^{r,k}$
satisfy the equations
\begin{subequations}
\label{eq:3:9:A_tau_prime_u_phi_z_eq_J_prime_u_phi_A_S_prime_u_phi_z_eq_J_prime_u_phi}
\begin{align}
\label{eq:3:9:A_tau_prime_u_phi_z_eq_J_prime_u_phi}
A_{\tau}^{\prime}(\primalu_{\tau})(\primalphi_{\tau},\dualz_{\tau})
& =
J^{\prime}(\primalu_{\tau})(\primalphi_{\tau})
\quad \hspace{0.3cm} \forall \primalphi_{\tau}\in X_{\tau}^{r}\,,
\\
\label{eq:3:9:A_S_prime_u_phi_z_eq_J_prime_u_phi}
A_{\tau h}^{\prime}(\primalu_{\tau h})(\primalphi_{\tau h},\dualz_{\tau h})
& =
J^{\prime}(\primalu_{\tau h})(\primalphi_{\tau h})
\quad \forall \primalphi_{\tau h}\in X_{\tau h}^{r,k}\,,
\end{align}
\end{subequations}
where $A_{\tau}^{\prime}(\cdot)(\cdot, \cdot)$ and
$A_{\tau h}^{\prime}(\cdot)(\cdot, \cdot)$ are given by
\begin{displaymath}
\begin{aligned}
A_{\tau}^{\prime}(\primalu_{\tau})(\primalphi_{\tau},\dualz_{\tau})  \coloneq
& \sum_{n=1}^N\int_{I_n}\big\{
(\primalpsi_{\tau},-\partial_{t} \dualvelocity_{\tau})
+ a^{\prime}(\primalu_{\tau})(\primalphi_{\tau},\dualz_{\tau})
\big\}\mathrm{d}t
\\
& - \sum_{n=1}^{N-1}(\primalpsi_{\tau,n}^-,[\dualvelocity_{\tau}]_{n})
+(\primalpsi_{\tau,N}^-,\dualvelocity_{\tau,N}^-)\,,
\end{aligned}
\end{displaymath}
and
\begin{displaymath}
\begin{aligned}
A_{\tau h}^{\prime}(\primalu_{\tau h})(\primalphi_{\tau h},\dualz_{\tau h}) \coloneq &
A_{\tau}^{\prime}(\primalu_{\tau h})(\primalphi_{\tau h},\dualz_{\tau h})
+{B}_{\Gamma_D}(\primalphi_{\tau h},\dualz_{\tau h})\,.
\end{aligned}
\end{displaymath}
\begin{remark}
We note that the directional derivatives of the Lagrangian functionals with
respect to their second argument lead to the primal problems given by
Eqs.~\eqref{eq:2:weakform},
\eqref{eq:8:time-discrete} and \eqref{eq:14:fully-discrete},
respectively.
\end{remark}
Now, using an abstract error representation result in terms of Lagrangian
functionals defined on different discretization levels given in \cite[Lemma~5.1]{BruchhaeuserBR12},
we get the following error representation formulas for the discretization errors
in space and time.
\begin{theorem}
\label{thm:1:ER}
Let $\{\primalu,\dualz\}\in X \times X$,
$\{\primalu_{\tau},\dualz_{\tau}\}
\in
X_{\tau}^{r} \times X_{\tau}^{r}$,
and
$\{\primalu_{\tau h},\dualz_{\tau h}\}
\in X_{\tau h}^{r,k} \times X_{\tau h}^{r,k}$
denote the stationary points of
$\mathcal{L}, \mathcal{L}_{\tau}$, and $\mathcal{L}_{\tau h}$
on the different levels of discretization, i.e.,
\begin{displaymath}
\begin{aligned}
\mathcal{L}^{\prime}(\primalu,\dualz)(\delta \primalu, \delta \dualz)
= \mathcal{L}_{\tau}^{\prime}(\primalu,\dualz)(\delta \primalu, \delta \dualz)
& = 0 \quad
\forall \{\delta \primalu,\delta \dualz\}\in X \times X\,,
\\
\mathcal{L}_{\tau}^{\prime}(\primalu_{\tau},\dualz_{\tau})
(\delta \primalu_{\tau}, \delta \dualz_{\tau})
& = 0
\quad \forall \{\delta \primalu_{\tau},\delta \dualz_{\tau}\}
\in X_{\tau}^{r} \times X_{\tau}^{r}\,,
\\
\mathcal{L}_{\tau h}^{\prime}(\primalu_{\tau h},\dualz_{\tau h})
(\delta \primalu_{\tau h}, \delta \dualz_{\tau h})
& = 0
\quad \forall \{\delta \primalu_{\tau h},\delta \dualz_{\tau h}\}
\in X_{\tau h}^{r,k} \times X_{\tau h}^{r,k}\,.
\end{aligned}
\end{displaymath}
Then, for the discretization errors in space and time we get the representation
formulas
\begin{subequations}
\label{eq:3:16}
\begin{align}
\label{eq:3:10a:J_u_minus_J_u_tau}
J(\primalu)-J(\primalu_{\tau}) = &
\frac{1}{2}\rho_{\tau}(\primalu_{\tau})(\dualz-\tilde{\dualz}_{\tau})
+ \frac{1}{2}\rho_{\tau}^{\ast}(\primalu_{\tau},\dualz_{\tau})
(\primalu-\tilde{\primalu}_{\tau})
+ \mathcal{R}_{\tau}\,,
\\
\label{eq:3:10b:J_u_tau_minus_J_u_tau_h}
J(\primalu_{\tau})-J(\primalu_{\tau h}) = &
\frac{1}{2}\rho_{\tau}(\primalu_{\tau h})(\dualz_{\tau}-\tilde{\dualz}_{\tau h})
+ \frac{1}{2}
\rho_{\tau}^{\ast}(\primalu_{\tau h},\dualz_{\tau h})
(\primalu_{\tau}-\tilde{\primalu}_{\tau h})
\\
\nonumber
&-\frac{1}{2}N_{\Gamma_D}(\velocity_{D,h}, \tilde{\dualz}_{\tau h}+\dualz_{\tau h})
\\
\nonumber
&
+ \frac{1}{2}{B}_{\Gamma_D}^{\prime}(\primalu_{\tau h}, \dualz_{\tau h})(\tilde{\primalu}_{\tau h}-\primalu_{\tau h},\tilde{\dualz}_{\tau h}-\dualz_{\tau h})
+{B}_{\Gamma_D}(\primalu_{\tau h}, \dualz_{\tau h})
+ \mathcal{R}_{h}\,,
\end{align}
\end{subequations}
where $\rho_{\tau}$ and $\rho_{\tau}^{\ast}$ denote the primal and dual residual
based on the semi-discrete schemes \eqref{eq:8:time-discrete} and \eqref{eq:3:9:A_tau_prime_u_phi_z_eq_J_prime_u_phi},
respectively, and ${B}_{\Gamma_D}^{\prime}(\cdot,\cdot)(\cdot,\cdot)$
denotes the G\^{a}teaux derivatives with respect to the first and second argument of the Nitsche terms ${B}_{\Gamma_D}(\cdot,\cdot)$ given in Eq.~\eqref{eq:13:nitsche}, defined as
\begin{subequations}
\label{eq:3:12:residuals}
\begin{align}
\label{eq:3:12:primal_residual}
\rho_{\tau}(\primalu)(\primalphi)\coloneq
\mathcal{L}^{\prime}_{\tau,\dualz}(\primalu,\dualz)(\primalphi)=&
F_{\tau}(\primalpsi)-A_{\tau}(\primalu)(\primalphi)\,,
\\
\label{eq:3:12:dual_residual}
\rho_{\tau}^{\ast}(\primalu,\dualz)(\primalphi)\coloneq
\mathcal{L}^{\prime}_{\tau,\primalu}(\primalu,\dualz)(\primalphi)=&
J^{\prime}(\primalu)(\primalphi)-A_{\tau}^{\prime}(\primalu)(\primalphi,\dualz)\,,
\\
\label{eq:Nitsche-Gateaux}
{B}_{\Gamma_D}^{\prime}(\primalu_{\tau h}, \dualz_{\tau h})(\tilde{\primalu}_{\tau h}-\primalu_{\tau h},\tilde{\dualz}_{\tau h}-\dualz_{\tau h})
\coloneq& {B}_{\Gamma_D}(\primalu_{\tau h},\tilde{\dualz}_{\tau h}-\dualz_{\tau h})
\\
\nonumber
&
+{B}_{\Gamma_D}(\tilde{\primalu}_{\tau h}-\primalu_{\tau h},\dualz_{\tau h})\,.
\end{align}
\end{subequations}
Here,
$\{\tilde{\primalu}_{\tau},\tilde{\dualz}_{\tau}\}\in X_{\tau}^{r}\times X_{\tau}^{r}$
and
$\{\tilde{\primalu}_{\tau h},\tilde{\dualz}_{\tau h}\} \in
X_{\tau h}^{r,k} \times X_{\tau h}^{r,k}$
can be chosen arbitrarily,
the discrete function $\velocity_{D,h}$ is an appropriate finite element approximation of an extension $\tilde{\velocity}_{D}$ in the sense that the trace equals $\velocity_{D}$ on $\Gamma_D$ in accordance with the trace theorem.
Finally, $\mathcal{R}_{\tau}$ and $\mathcal{R}_{h}$ are
higher-order remainder terms with respect to the errors $\primalu-\primalu_\tau,\dualz-\dualz_\tau$ and
$\primalu_{\tau}-\primalu_{\tau h},\dualz_{\tau}-\dualz_{\tau h}$, cf.~\cite{BruchhaeuserBR12,BruchhaeuserBR01}.
\end{theorem}
\begin{mproof}
The technique to prove the temporal and spatial error representation formulas
goes back to an abstract error result in terms of Lagrangian functionals defined
on different discretization levels given by Besier and Rannacher applied to the
incompressible Navier--Stokes equations using local projection stabilization (LPS)
in\cite{BruchhaeuserBR12}.
More precisely, we are using the abstract result given in \cite[Lemma~5.1]{BruchhaeuserBR12}
with the following settings:
\begin{displaymath}
\begin{array}{l@{\,}c@{\,}l@{\,}l@{\,}c@{\,}l@{\,}l@{\,}c@{\,}l@{\,}l@{\,}c@{\,}l@{\,}l@{\,}c@{\,}l@{\,}}
L & = & \mathcal{L}\,,\;\; &
\tilde{L}&=&\mathcal{L}_{\tau}\,,\;\;&
Y_1 &=& X \times X\,,\;\;&
Y_2 &=& X_{\tau}^{r} \times X_{\tau}^{r}\,,\;\;&
Y&\coloneq&Y_1+Y_2\,,
\\
L & = & \mathcal{L}_{\tau}\,,\;\; &
\tilde{L}&=&\mathcal{L}_{\tau h}\,,\;\;&
Y_1 &=& X_{\tau}^{r} \times X_{\tau}^{r}\,,\;\;&
Y_2 &=& X_{\tau h}^{r,k} \times X_{\tau h}^{r,k}\,,\;\;&
Y&\coloneq&Y_1\,,
\end{array}
\end{displaymath}
where $\mathcal{L},\mathcal{L}_{\tau},$ and $\mathcal{L}_{\tau h}$ are the
Lagrangian functionals defined in
(\ref{eq:3:3:Lagrangians}) and $Y,Y_1$ and $Y_2$
are function spaces defined in \cite[Lemma~5.1]{BruchhaeuserBR12}.
The resulting Nitsche terms are treated in the same way as the stabilization terms in the proof of \cite[Thm.~5.2]{BruchhaeuserBR12}.
\end{mproof}


\section{Practical Aspects}
\label{sec:4:practical}
In this section, we introduce local error indicators in space and time, present
the space-time adaptive algorithm, and describe the geometric multigrid
preconditioned Newton--Krylov solver together with selected aspects of the
parallel implementation.

\subsection{Local Error Indicators in Space and Time}

For the marking step in the adaptive mesh refinement process, we require local
versions of the error representations \eqref{eq:3:16} from
Theorem~\ref{thm:1:ER}. Neglecting the higher-order remainder terms
\(\mathcal R_\tau\) and \(\mathcal R_h\) in \eqref{eq:3:16} yields the global a
posteriori estimators
\begin{align*}
J(\primalu)-J(\primalu_{\tau})&\approx \eta_{\tau}\,,
\\
J(\primalu_{\tau})-J(\primalu_{\tau h})&\approx \eta_{h}\,,
\end{align*}
with
\begin{align*}
\eta_{\tau}\coloneq& \frac{1}{2}\rho_{\tau}(\primalu_{\tau})(\dualz-\tilde{\dualz}_{\tau})
+ \frac{1}{2}\rho_{\tau}^{\ast}(\primalu_{\tau},\dualz_{\tau})
(\primalu-\tilde{\primalu}_{\tau})\,,
\\[1ex]
\eta_{h}\coloneq& \frac{1}{2}\rho_{\tau}(\primalu_{\tau h})(\dualz_{\tau}-\tilde{\dualz}_{\tau h})
+ \frac{1}{2}
\rho_{\tau}^{\ast}(\primalu_{\tau h},\dualz_{\tau h})
(\primalu_{\tau}-\tilde{\primalu}_{\tau h})
\\
&
-\frac{1}{2}N_{\Gamma_D}(\velocity_{D,h}, \tilde{\dualz}_{\tau h}+\dualz_{\tau h})
+ \frac{1}{2}{B}_{\Gamma_D}^{\prime}(\primalu_{\tau h}, \dualz_{\tau h})(\tilde{\primalu}_{\tau h}-\primalu_{\tau h},\tilde{\dualz}_{\tau h}-\dualz_{\tau h})
\\
&+{B}_{\Gamma_D}(\primalu_{\tau h}, \dualz_{\tau h})
\,.
\end{align*}

To obtain implementable estimators, we replace the unknown exact and
semi-discrete solutions by the fully discrete approximations
\(\primalu_{\tau h},\dualz_{\tau h}\in X_{\tau h}^{r,k}\). The temporal weights
\(\primalu-\tilde{\primalu}_{\tau}\) and
\(\dualz-\tilde{\dualz}_{\tau}\), and the spatial weights
\(\primalu_{\tau}-\tilde{\primalu}_{\tau h}\) and
\(\dualz_{\tau}-\tilde{\dualz}_{\tau h}\), are then approximated as follows.

\paragraph{Temporal weights.}
We approximate the temporal weights by higher-order reconstructions in time:
\begin{displaymath}
\primalu-\tilde{\primalu}_{\tau} \approx \mathrm{E}_{\tau}^{(r+1)}\primalu_{\tau h}-\primalu_{\tau h}\,,
\qquad
\dualz-\tilde{\dualz}_{\tau} \approx \mathrm{E}_{\tau}^{(r+1)}\dualz_{\tau h}-\dualz_{\tau h}\,.
\end{displaymath}
Here \(\mathrm{E}_{\tau}^{(r+1)}\) denotes a reconstruction operator in time
that lifts the discrete solution to a piecewise polynomial of degree \(r+1\),
see \cite{BruchhaeuserBB24,BruchhaeuserB22}.

\paragraph{Spatial weights.}
We approximate the spatial weights by patch-wise higher-order interpolation:
\begin{equation}
  \primalu_{\tau}-\tilde{\primalu}_{\tau h}  \approx  \mathrm{I}_{2h}^{(2k)}\primalu_{\tau h}-\primalu_{\tau h}\,,
  \qquad
  \dualz_{\tau}-\tilde{\dualz}_{\tau h}  \approx  \mathrm{I}_{2h}^{(2k)}\dualz_{\tau h}-\dualz_{\tau h}\,.
\end{equation}
Here \(\mathrm{I}_{2h}^{(2k)}\) denotes a patch-wise interpolation operator that
acts on a patch of size \(2h\) and lifts the discrete solution to a piecewise
polynomial of degree \(2k\), see \cite{BruchhaeuserBBEMTW25}.

With these approximations, the implemented temporal and spatial estimators are
given by
\begin{subequations}\label{eq: implemented error estimator}
  \begin{align}\label{eq: implemented error estimator-tau}
\tilde{\eta}_{\tau}\approx&\frac{1}{2}\rho_{\tau}(\primalu_{\tau h})(\mathrm{E}_{\tau}^{(r+1)}\dualz_{\tau h}-\dualz_{\tau h})
+ \frac{1}{2}\rho_{\tau}^{\ast}(\primalu_{\tau h},\dualz_{\tau h})(\mathrm{E}_{\tau}^{(r+1)}\primalu_{\tau h}-\primalu_{\tau h})\,,
\\
\label{eq: implemented error estimator-h}
\tilde{\eta}_{h}\approx&\frac{1}{2}\rho_{\tau}(\primalu_{\tau h})(\mathrm{I}_{2h}^{(2k)}\dualz_{\tau h}-\dualz_{\tau h})
+ \frac{1}{2}\rho_{\tau}^{\ast}(\primalu_{\tau h},\dualz_{\tau h})(\mathrm{I}_{2h}^{(2k)}\primalu_{\tau h}-\primalu_{\tau h})
\\\nonumber
&
-\frac{1}{2}N_{\Gamma_D}(\velocity_{D,h}, \mathrm{I}_{2h}^{(2k)}\dualz_{\tau h}+\dualz_{\tau h})
+ \frac{1}{2}{B}_{\Gamma_D}^{\prime}(\primalu_{\tau h}, \dualz_{\tau h})(\mathrm{I}_{2h}^{(2k)}\primalu_{\tau h}-\primalu_{\tau h},\mathrm{I}_{2h}^{(2k)}\dualz_{\tau h}-\dualz_{\tau h})
\\
\nonumber
&
+{B}_{\Gamma_D}(\primalu_{\tau h}, \dualz_{\tau h})
    \,.
  \end{align}
\end{subequations}

For mesh marking, these implemented estimators are decomposed into local contributions applying integration by parts on each element $K\in\mathcal{T}_h^n$ to the viscous term.
Following \cite{BruchhaeuserBR01,BruchhaeuserBR03,BruchhaeuserBR12}, we write
\begin{equation}
  \label{eq:localized-errors-iso}
\tilde{\eta}_{\tau}=\sum_{n=1}^{N}\sum_{K\in\mathcal{T}_h^n}\tilde{\eta}_{\tau}^{K,n},
\qquad
\tilde{\eta}_h=\sum_{n=1}^{N}\sum_{K\in\mathcal{T}_h^n}\tilde{\eta}_{h}^{K,n}\,.
\end{equation}
The local indicators \(\tilde{\eta}_{\tau}^{K,n}\) and \(\tilde{\eta}_{h}^{K,n}\) are then used
for slab and element marking.

\subsection{Space-Time Adaptive Algorithm}
\label{sec:4:2:algorithm}

We present a condensed version of the adaptive solution algorithm based on
tensor-product space-time slabs \(Q_n\). For further details on the slab data
structures, we refer to \cite{BruchhaeuserBKB22} and
\cite[Sec.~4.4]{BruchhaeuserB22}.

\begin{algorithm}[htb]
  \caption{Goal-Oriented Space-Time Adaptivity}
  \label{alg:goal-adaptivity}
  \begin{algorithmic}[1]
    \State \textbf{Input:} Initial space-time mesh \(Q_n^1=\mathcal{T}_{h,n}^1\times I_n^1\), \(n=1,\dots,N^1\), and parameter \(\omega\in[1.5,3.5]\).
    \State Set iteration counter \(\iota \gets 1\).
    \Repeat
      \State Solve the \textbf{primal} problem \eqref{eq:14:fully-discrete} to obtain \(\boldsymbol u_{\tau h}^{\iota}\).
      \State Solve the \textbf{dual} problem \eqref{eq:3:9:A_S_prime_u_phi_z_eq_J_prime_u_phi} to obtain \(\boldsymbol z_{\tau h}^{\iota}\).
      \State Evaluate the global space-time error indicators \(\tilde{\eta}_{\tau}^{\iota}\) and \(\tilde{\eta}_{h}^{\iota}\).
      \If{\(|\tilde{\eta}_{\tau}^{\iota}|>\omega\,|\tilde{\eta}_{h}^{\iota}|\)}
        \State Adapt the \textbf{temporal} mesh only.
      \ElsIf{\(|\tilde{\eta}_{h}^{\iota}|>\omega\,|\tilde{\eta}_{\tau}^{\iota}|\)}
        \State Adapt the \textbf{spatial} mesh only.
      \Else
        \State Adapt both the \textbf{temporal} and \textbf{spatial} meshes.
      \EndIf
      \State Set \(\iota \gets \iota+1\).
    \Until{a stopping criterion is met.}
    \State \textbf{Output:} Adapted space-time mesh and final primal solution.
  \end{algorithmic}
\end{algorithm}

\begin{remark}[Adaptation of the space-time mesh]\label{rem:adaptation}
  In Algorithm~\ref{alg:goal-adaptivity}, we do not specify the concrete
  marking strategy. We briefly describe the strategy used in this work.

  Let \(\omega\in[1.5,3.5]\) be the criterion for temporal or spatial error
  domination. Let \(\theta_{\mathrm{top}}^\tau\in(0,1]\) denote the fraction of
  slabs to refine in time. Let
  \(\theta_{\mathrm{top},1}^{h,n},\theta_{\mathrm{top},2}^{h,n}\in(0,1]\) be
  the fractions of elements to refine on slab \(n\) without and with temporal
  refinement, respectively. Let
  \(\theta_{\mathrm{bot},1}^{h},\theta_{\mathrm{bot},2}^{h}\in[0,1)\) be the
  corresponding fractions of elements to coarsen.

  Let \(\tilde\eta_\tau\) and \(\tilde\eta_h\) denote the global
  temporal and spatial indicators used for marking. The marking strategy is as
  follows.
  \begin{itemize}
    \item If \(|\tilde{\eta}_\tau|\ge \omega\,|\tilde{\eta}_h|\), then we
    mark the top \(\theta_{\mathrm{top}}^\tau\) fraction of slabs for temporal
    refinement.
    \item If \(|\tilde{\eta}_h|\ge \omega\,|\tilde{\eta}_\tau|\), then on
    each slab \(n\),
    \begin{itemize}
      \item we refine the top \(\theta_{\mathrm{top},1}^{h,n}\) fraction of
      elements if the slab is not marked for temporal refinement,
      \item we refine the top \(\theta_{\mathrm{top},2}^{h,n}\) fraction of
      elements if the slab is marked for temporal refinement,
      \item we coarsen the bottom \(\theta_{\mathrm{bot},1}^{h}\) fraction of
      elements if the slab is not marked for temporal refinement,
      \item we coarsen the bottom \(\theta_{\mathrm{bot},2}^{h}\) fraction of
      elements if the slab is marked for temporal refinement.
    \end{itemize}
    \item If neither indicator dominates, then we apply both the temporal and
    spatial marking rules above.
  \end{itemize}
\end{remark}

In Remark~\ref{rem:adaptation}, we introduced slab-wise refinement fractions. In
this work, we either keep these fractions constant, as in previous work, or
adapt them according to the following rule.

\begin{remark}[Slab-wise adaptation of the spatial refinement fractions]
  On iteration \(\iota\), define on each space-time slab the local indicator
  \[
    \tilde{\eta}_{\tau h}^{n,\iota}=\tilde{\eta}_{h}^{n,\iota}+\tilde{\eta}_{\tau}^{n,\iota}\,,
  \]
  and the corresponding global quantity
  \[
    \tilde{\eta}_{\tau h}^{\iota}=\sum_{n=1}^{N^\iota}\tilde{\eta}_{\tau h}^{n,\iota}\,.
  \]
  If \(\tilde{\eta}_{\tau h}^{\iota}=0\), we keep the nominal refinement fractions.
  Otherwise, we define the normalized slab-wise weights
  \[
    \alpha_n=\frac{\tilde{\eta}_{\tau h}^{n,\iota}}{\tilde{\eta}_{\tau h}^{\iota}}\,,
    \qquad
    \alpha_{\min}=\min_n\alpha_n\,,
    \qquad
    \alpha_{\max}=\max_n\alpha_n\,.
  \]
  These weights are used to adapt the nominal fraction
  \(\theta_{\mathrm{top},2}^{h}\) on slabs that are marked both in time and in
  space. To avoid large deviations from the nominal value, we bound the adapted
  fractions by
  \[
    \theta_{\mathrm{top},2}^{h,-}=\tfrac34\,\theta_{\mathrm{top},2}^{h}\,,
    \qquad
    \theta_{\mathrm{top},2}^{h,+}=\min\!\left\{1,\tfrac43\,\theta_{\mathrm{top},2}^{h}\right\}\,.
  \]
  The slab-wise fraction \(\theta_{\mathrm{top},2}^{h,n}\) is then defined by
  \[
    \theta_{\mathrm{top},2}^{h,n}=
    \begin{cases}
      \theta_{\mathrm{top},2}^{h}, & \alpha_{\max}-\alpha_{\min}<\varepsilon,\\[1ex]
      \displaystyle
      \theta_{\mathrm{top},2}^{h,-}
      +\bigl(\theta_{\mathrm{top},2}^{h,+}-\theta_{\mathrm{top},2}^{h,-}\bigr)
      \frac{\alpha_n-\alpha_{\min}}{\alpha_{\max}-\alpha_{\min}},
      & \text{otherwise},
    \end{cases}
  \]
  where \(\varepsilon>0\) is a small tolerance. Hence
  \[
    \theta_{\mathrm{top},2}^{h,n}\in
    \bigl[\theta_{\mathrm{top},2}^{h,-},\,\theta_{\mathrm{top},2}^{h,+}\bigr]\,,
  \]
  and the spatial refinement fraction is distributed according to the relative
  slab-wise error contributions.
\end{remark}

\begin{remark}[Dynamic meshes]
We solve the fully problem \eqref{eq:14:fully-discrete} on the mesh slabs $Q_n= \mathcal T_h^n \times I_n$. Here, the mesh $\mathcal T_h^n$ and the pair of finite element spaces $
\boldsymbol V_h^k(\mathcal T_h^n)\times Q_h^{k-1}(\mathcal T_h^n)$
are allowed to change from subinterval to subinterval. The velocity field $\velocity_{n-1}^{-}\coloneq \velocity_{\tau,h}(t_{n-1})\in \boldsymbol V_h^k(\mathcal T_h^{n-1})$, computed on $I_{n-1}$, then satisfies the  discrete divergence-free constraint with respect to all test functions in $Q_h^{k-1}(\mathcal T_h^{n-1})$ and violates the constraint with respect to $Q_h^{k-1}(\mathcal T_h^{n})$ on $I_n$, if \(\mathcal T_h^{n-1}\neq \mathcal T_h^n\). To recover the discrete divergence-free constraint of the coupling term on $I_n$, we substitute $\velocity_{n-1}^{-}$ by its Helmholtz projection $\Pi_h^n\velocity_{n-1}^{-}\in \boldsymbol V_h^k(\mathcal T_h^{n})$, satisfying with some additional pressure function $\vartheta_h^n\in Q_h^{k-1}(\mathcal T_h^n)$ the variational equation \cite{besier_pressure_2012,10.1093/imanum/drt001},
\begin{equation}
\label{eq-stokes-projection}
(\Pi_h^n\velocity_{n-1}^{-},\boldsymbol v_h)
-(\vartheta_h^n,\nabla\cdot \boldsymbol v_h)
+(\nabla\cdot \Pi_h^n\velocity_{n-1}^{-},q_h)
=
(\widehat{\velocity}_{\tau h,n-1}^{\,n},\boldsymbol v_h)
\end{equation}
for all
\((\boldsymbol v_h,q_h)\in \boldsymbol V_h^k(\mathcal T_h^n)\times
Q_h^{k-1}(\mathcal T_h^n)\). The projection $\Pi_h^n\velocity_{n-1}^{-}\in \boldsymbol V_h^k(\mathcal T_h^{n})$ is discretely divergence free with respect to the decomposition $\mathcal T_h^n$ of $\Omega$,
\[
(\nabla\cdot \Pi_h^n\velocity_{n-1}^{-},q_h)=0
\]
for all $q_h\in Q_h^{k-1}(\mathcal T_h^n)$. Alternatively, a Stokes projection can be used. As studied in \cite{10.1093/imanum/drt001}, either projections can lead to the accumulation of projection errors that are due to \eqref{eq-stokes-projection}. Instabilities and failure of can occur if arbitrary mesh modifications without hierarchical embeddings along with a high number of changes are involved. In our numerical computations presented in Section~\ref{sec:5:numerical}, we do not observe such shortcomings. Suitable space-time projections or further post-processing might overcome accumulating projection errors, even for high numbers of mesh changes.
\end{remark}

\subsection{Geometric Multigrid Preconditioned Newton-Krylov Method}

\begin{algorithm}[htb]
  \caption{Geometric multigrid method for the linear system
  \(A_Lx_L=b_L\) on the finest level \(L\). Under the standard smoothing and
  approximation properties, a multigrid V-cycle can achieve \(O(N)\)
  complexity.}
  \label{alg:gmg}
  \begin{algorithmic}[1]
    \Procedure{multigrid}{$\ell,\boldsymbol A_\ell,\boldsymbol b_\ell,\boldsymbol x_\ell$}
    \State \(\boldsymbol s_\ell \leftarrow \mathcal S(\boldsymbol A_\ell,\boldsymbol b_\ell,\boldsymbol x_\ell)\)\Comment{Pre-smoothing}
    \State \(\boldsymbol r_\ell \leftarrow \boldsymbol b_\ell - \boldsymbol A_\ell \boldsymbol s_\ell\)\Comment{Residual}
    \State \(\boldsymbol r_{\ell-1} \leftarrow \mathcal R(\ell,\boldsymbol r_\ell)\)\Comment{Restriction}
    \If{\(\ell=1\)}
      \State \(\boldsymbol c_0 \leftarrow \boldsymbol A_0^{-1}\boldsymbol r_0\)\Comment{Coarse-grid solve}
    \Else
      \State \(\boldsymbol c_{\ell-1} \leftarrow \Call{multigrid}{\ell-1,\boldsymbol A_{\ell-1},\boldsymbol r_{\ell-1},\boldsymbol 0}\)
    \EndIf
    \State \(\boldsymbol x'_\ell \leftarrow \boldsymbol s_\ell + \mathcal P(\ell,\boldsymbol c_{\ell-1})\)\Comment{Prolongation}
    \State \(\boldsymbol s'_\ell \leftarrow \mathcal S(\boldsymbol A_\ell,\boldsymbol b_\ell,\boldsymbol x'_\ell)\)\Comment{Post-smoothing}
    \State \Return \(\boldsymbol s'_\ell\)
    \EndProcedure
  \end{algorithmic}
\end{algorithm}

\begin{algorithm}[htb]
  \caption{Newton--GMRES method with GMG preconditioning}
  \label{alg:newton-gmres}
  \begin{algorithmic}[1]
    \State \textbf{Input:} initial guess \(\boldsymbol x^{(0)}\), tolerance \(\epsilon\), maximum number of iterations \(N_{\max}\), reassembly threshold \(\rho_{\max}\).
    \State Assemble the initial Jacobian \(J\gets \mathcal A'(\boldsymbol x^{(0)})\).
    \For{\(l=0,1,\dots,N_{\max}-1\)}
      \State Compute the residual \(\boldsymbol r^{(l)}\gets \boldsymbol F_{\tau h}-\boldsymbol{\mathcal A}(\boldsymbol x^{(l)})\).
      \If{\(\|\boldsymbol r^{(l)}\|<\epsilon\)}
        \State \Return \(\boldsymbol x^{(l)}\)
      \EndIf
      \If{\(l>0\) \textbf{and} \(\|\boldsymbol r^{(l)}\|/\|\boldsymbol r^{(l-1)}\|>\rho_{\max}\)}
        \State Reassemble \(J\gets \mathcal A'(\boldsymbol x^{(l)})\).
      \EndIf
      \State Solve \(J\,\boldsymbol w^{(l)}=\boldsymbol r^{(l)}\) by GMRES preconditioned by GMG, see Algorithm~\ref{alg:gmg}.
      \State Update \(\boldsymbol x^{(l+1)}\gets \boldsymbol x^{(l)}+\boldsymbol w^{(l)}\).
    \EndFor
    \State \textbf{Output:} \(\boldsymbol x^{(N_{\max})}\).
  \end{algorithmic}
\end{algorithm}

The space-time finite element discretization of the incompressible
Navier--Stokes equations leads to the nonlinear algebraic
system \eqref{eq:14:fully-discrete}.

To solve this problem, we employ a Newton--Krylov method preconditioned by
geometric multigrid. Newton's method is applied in the weak form
\begin{equation}\label{newton}
  A_{\tau h}'(\boldsymbol u^{(l-1)})(\boldsymbol w^{(l)},\boldsymbol\varphi_{\tau h})
  =
  F_{\tau h}(\boldsymbol\varphi_{\tau h})
  -
  A_{\tau h}(\boldsymbol u^{(l-1)})(\boldsymbol\varphi_{\tau h}),
  \quad
  \boldsymbol u^{(l)}=\boldsymbol u^{(l-1)}+\boldsymbol w^{(l)}.
\end{equation}
Here \(A_{\tau h}'(\boldsymbol u^{(l-1)})\) denotes the Fr\'echet derivative of
the semilinear form \(A_{\tau h}\) at \(\boldsymbol u^{(l-1)}\). Algebraically,
this corresponds to the Jacobian matrix evaluated at
\(\boldsymbol u^{(l-1)}\) and applied to the Newton correction
\(\boldsymbol w^{(l)}\). The Jacobian can be computed analytically, see
\cite[Sec.~4.4.2]{BruchhaeuserR17}. To reduce numerical costs, the Jacobian is reused over
several Newton iterations and consecutive slabs, and is reassembled only when
the observed residual reduction factor \(\|\boldsymbol r^{(l)}\|/\|\boldsymbol r^{(l-1)}\|\) exceeds a prescribed threshold \(\rho_{\max}\). In our computations,
\(\rho_{\max}\) is chosen empirically and typically lies between \(0.05\) and
\(0.1\).

Each Newton step requires the solution of a large, sparse, nonsymmetric, and
indefinite linear system arising from the saddle-point structure of the
Navier--Stokes equations. We use the Generalized Minimal Residual method
(GMRES) \cite{KelleyIterativeMethodsLinear1995,Saad1996} to solve these linear
systems. To accelerate convergence, GMRES is preconditioned by geometric
multigrid. On each slab $I_n$, the multigrid method is built on a hierarchy of nested
discrete saddle-point spaces
\[
  \mathcal X_0^n \subset \mathcal X_1^n \subset \cdots \subset \mathcal X_L^n,
\]
associated with nested spatial meshes
\(\mathcal T_0^n,\dots,\mathcal T_L^n=\mathcal T_h^n\). On each level,
high-frequency errors are damped by a smoothing operator, while low-frequency
components are transferred to coarser levels by \(L^2\)-projection-based
restriction and interpolative prolongation. On the coarsest level, a direct
solver is employed. Under the standard multigrid assumptions, one V-cycle has
linear complexity in the number of degrees of freedom.

To treat the saddle-point structure efficiently and to obtain a parallelizable
smoother, we employ a Vanka-type smoother \cite{Vanka1985} within the multigrid
cycle. It is defined by
\begin{equation}\label{vanka0}
  \mathcal{V}_\ell(A_\ell)
  =
  \sum_{K\in\mathcal T_\ell}
  R_K^\top \bigl[R_KA_\ell R_K^\top\bigr]^{-1}R_K,
\end{equation}
where \(A_\ell\) denotes the system matrix on level \(\ell\), and \(R_K\) is
the restriction operator to the degrees of freedom associated with the element
\(K\in\mathcal T_\ell\). One smoothing step is then given by
\begin{equation}\label{vanka}
  \mathcal S(\boldsymbol A_\ell,\boldsymbol b_\ell,\boldsymbol x_k)
  =
  \boldsymbol x_k
  +
  \omega\,\mathcal V_\ell(A_\ell)\,
  (\boldsymbol b_\ell-\boldsymbol A_\ell\boldsymbol x_k),
\end{equation}
with damping parameter \(\omega\in(0,1]\). In our experience, two to five pre-
and post-smoothing steps are sufficient for robust convergence. The
discontinuous pressure approximation is beneficial for the Vanka smoother, since
it localizes the pressure degrees of freedom and thereby improves the
efficiency of the element-wise local solves.

\begin{remark}[Weighting in the Vanka smoother]
  In the implementation, the additive operator \eqref{vanka0} is combined with
  a diagonal partition-of-unity weighting. We choose the weight of each degree
  of freedom as the reciprocal of its valence. An alternative is a restricted
  Vanka smoother, which assigns each degree of freedom to a unique element. This
  reduces communication, but in our tests it led to a less effective smoother
  and hence to a larger number of GMRES iterations. Since in the present
  matrix-based setting the dominant cost comes from global matrix-vector
  products, we use the weighted variant here.
\end{remark}

\begin{remark}[Dual linear systems]
For fixed primal solution, the dual problem is linear and solved backward in
time. Its slab-wise discretization leads to nonsymmetric saddle-point systems
with diffusion, backward transport, and a reaction term. In the
implementation, these systems are treated with the same GMRES-preconditioned
geometric multigrid method as the primal Newton linearizations, based on the
same level hierarchy, intergrid transfers, and Vanka-type smoothing.
\end{remark}

In summary, the nonlinear space-time finite element system is linearized by
Newton's method and the resulting linear systems are solved by GMRES
preconditioned by geometric multigrid with a Vanka-type smoother. For a more
detailed description of the solver architecture and its parallelization, we
refer to \cite{BruchhaeuserAB23}.

\subsection{Parallel implementation}

The implementation of the adaptive algorithm is distributed by MPI and uses the
parallel infrastructure of \texttt{deal.II} for error estimation, marking,
refinement, and multigrid setup across all processes. No stage is executed on a
dedicated master rank. Collective communication is still required for the
assembly of global indicators and for the coarse-grid solve.

\begin{remark}[MPI-distributed adaptation]
\leavevmode
\begin{itemize}
\item \emph{Local error estimation with face terms.} Each rank computes its
  element-wise residual contributions, including jump terms on interior faces.
  For faces on process boundaries, ghost-element values are exchanged with
  neighboring ranks so that these face terms are evaluated consistently. Each
  face contribution is computed exactly once and communicated if necessary.

\item \emph{Global indicator assembly.} The local error contributions per slab
  or per element are summed by \texttt{MPI\_Allreduce} to form the global
  indicators \(\tilde{\eta}_h\) and \(\tilde{\eta}_\tau\), and to normalize
  the slab-wise weights \(\alpha_n\).

\item \emph{Marking for refinement and coarsening.} Each rank marks its locally
  owned slabs or elements according to the global decision between
  \(\tilde{\eta}_h\) and \(\tilde{\eta}_\tau\) and according to the local
  fractions \(\theta_{\mathrm{top},2}^{h,n}\). No additional communication is
  required for the marking itself.

\item \emph{Distributed mesh refinement.} At present, we do not repartition the
  mesh, so load imbalance may accumulate over several adaptive cycles.
  Repartitioning is available in \texttt{deal.II}. In the present space-time
  setting, however, the many mesh transfers make it less clear how much
  performance can be gained. The development of an optimized repartitioning
  strategy is left to future work.

\item \emph{Multigrid preconditioner with local smoothing.} On each rank, the
  space-time element-wise Vanka smoother is applied to the locally owned elements on
  each level, with the usual ghost-layer exchange where required. The
  coarse-grid correction is handled by a global, direct solve.
\end{itemize}
\end{remark}

At present, load balancing and repartitioning are not used. To maintain
scalability after several adaptive cycles, we plan to incorporate dynamic
repartitioning in future work.

\section{Numerical Examples}
\label{sec:5:numerical}
In this section, we investigate the performance properties of the proposed approach, which integrates goal-oriented error control with iterative solver techniques. In particular, we study the convergence behavior, computational efficiency, as well as stability properties of the underlying space-time adaptive algorithm (Algorithm~\ref{alg:goal-adaptivity}) introduced in Sect.~\ref{sec:4:2:algorithm}.
These properties are examined by means of the following three numerical examples.
\begin{enumerate}
\item \emph{Exact Analytical Solution:}
The first example is an academic test problem including a given analytical solution that goes back to \cite[Ex.~1]{BruchhaeuserBR12}.
It is used to validate accuracy and efficiency of the derived approach for different discretizations up to polynomial order 3 in time and 4 in space.
\item \emph{DFG Benchmark 2D-3:} The second example is the famous DFG flow around cylinder benchmark 2D-3 that simulates the time-periodic behavior of a fluid in a pipe with a circular obstacle that goes back to \cite{BruchhaeuserST96}.
The goal functional is chosen to control the mean drag.
\end{enumerate}
%
We evaluate the accuracy and efficiency of our proposed approach in terms of several quantities. Firstly, to assess the accuracy of the duality based space-time error estimator, we compute the effectivity index defined as the ratio of the estimated to the exact error, given by
\begin{equation}
\label{eq:Ieff}
\Ieff\coloneq \left|\frac{
\tilde{\eta}_{\tau}
+
\tilde{\eta}_{h}}
{J(\primalu)-J(\primalu_{\tau h})}\right|\,.
\end{equation}
%
Secondly, to measure the gain in accuracy over the increase of work, we determine the empirical convergence rates, averaged over all DWR-loops,
\begin{equation}\label{eq:mean-rates}
  \overline{\operatorname{rate}}_{N_{\mathrm{ST}}}^{\|e_T\|}
  \quad \text{and} \quad
  \overline{\operatorname{rate}}_{W_{\mathrm{tot}}}^{\|e_T\|}
\end{equation}
of the reduction of the error in the goal quantity with respect to the increase of the number of space-time degrees of freedom \(N_{\mathrm{ST}}\) and the total work \(W_{\mathrm{tot}}\). The latter is measured by the cumulative linear solver work
\begin{equation}\label{eq:work}
  W_{\mathrm{tot}} = W_p + W_d,\quad\text{with primal and dual contributions}
\end{equation}
\begin{equation}
   W_p \coloneq \sum_{n=1}^{N_{\mathrm{time}}} n_{\mathrm{L},n}^{p}  N_{\mathrm{DoF},n}^{p},
  \;
  W_d \coloneq \sum_{n=1}^{N_{\mathrm{time}}} n_{\mathrm{L},n}^{d}  N_{\mathrm{DoF},n}^{d}\,,
\end{equation}
where \(n_{\mathrm{L},n}^{p}\) and \(n_{\mathrm{L},n}^{d}\) denote the numbers
of primal and dual linear iterations on slab \(Q_n = \mathcal{T}_h^n \times I_n\), and
\(N_{\mathrm{DoF},n}^{p}\), \(N_{\mathrm{DoF},n}^{d}\) the corresponding numbers
of spatial degrees of freedom. The nonlinear primal problem is solved by Newton's method. We therefore also report the number of nonlinear iterations on slab \(Q_n\), denoted by \(n_{\mathrm{NL},n}\). Accordingly, we let \(n_{\mathrm{L},n}^{p}\) denote the total number of linear iterations accumulated over the \(n_{\mathrm{NL},n}\) Newton steps on slab \(Q_n\).
By \(n_{\mathrm{NL}}\), \(n_{\mathrm{L}}^{p}\), and
\(n_{\mathrm{L}}^{d}\) we denote the total numbers of nonlinear, primal linear, and dual linear iterations over all slabs within one DWR loop $\ell$. In addition, we record
\begin{equation}\label{eq:iteration-indicators}
  \min_{\ell=1,\dots,\ell_{\max}} \Bigl(\frac{n_{\mathrm{NL}}}{N_{\mathrm{time}}}\Bigr),\;\;
  \max_{\ell=1,\dots,\ell_{\max}} \Bigl(\frac{n_{\mathrm{NL}}}{N_{\mathrm{time}}}\Bigr),\;\;
  \min_{\ell=1,\dots,\ell_{\max}} \Bigl(\frac{n_{\mathrm{L}}^{p}}{n_{\mathrm{NL}}}\Bigr),\;\;
  \max_{\ell=1,\dots,\ell_{\max}} \Bigl(\frac{n_{\mathrm{L}}^{p}}{n_{\mathrm{NL}}}\Bigr),
\end{equation}
together with the relative dual contribution
\begin{equation}\label{eq:dual-work-fraction}
  \frac{W_d}{W_{\mathrm{tot}}}.
\end{equation}
Throughout, by $\ell$ we denote the refinement level or DWR loop number, $N_{\text{time}}$ the number of subintervals in time, $N_{\text{space}}^{\text{max}}$ the maximum number of degrees of freedom associated with a spatial triangulation within one loop and $N_{\text{tot}}$
the number of total degrees of freedom in space and time.
We consider different finite element approximations dG($r$)-$Q_k/P_{k-1}^{\mathrm{disc}}$ denoting a discontinuous in time approximation of polynomial degree $r\geq 0$ as well as a continuous velocity and discontinuous pressure approximations in space of polynomial degree $k\geq2$.

The implementation is based on the \texttt{deal.II} finite element library~\cite{BruchhaeuserABBFGHHKKMMPPSTWZ24}. The tests are run on a single node with 2 Intel Xeon Platinum 8360Y CPUs and \SI[scientific-notation=false,round-precision=0]{1024}{\giga\byte} RAM of the HPC cluster HSUper at Helmut Schmidt University.

\subsection{Example 1: Exact Analytical Solution} %
As a first test case, we consider a model problem on the space-time domain $\Omega\times I = [0,1]^2\times [0, 1]$ with prescribed solution given for the velocity $\boldsymbol v \colon \Omega\times I \to \mathbb{R}^2$ and pressure
$p \colon \Omega\times I \to \mathbb{R}$ by
\begin{equation}
\label{eq:conv-test}
\begin{aligned}
\boldsymbol v(\boldsymbol{x},\,t) &= \sin(t) \begin{pmatrix} \sin^2(\pi x) \sin(\pi y) \cos(\pi y) \\
\sin(\pi x) \cos(\pi x) \sin^2(\pi y) \end{pmatrix},
\\
p(\boldsymbol{x},\,t)
&= \sin(t) \sin(\pi x) \cos(\pi x)
  \sin(\pi y) \cos(\pi y)\,.
\end{aligned}
\end{equation}
We choose the kinematic viscosity as $\nu \in \{1, 10^{-3}\}$ and choose
the external force $\boldsymbol{f}$ such that the solution~\eqref{eq:conv-test} satisfies~\eqref{eq:1:navier}. The initial velocity is prescribed as zero and homogeneous Dirichlet boundary conditions are imposed on
$\partial\Omega$ for all times
\[
  \boldsymbol v = \boldsymbol{0}\text{ on }\Omega\times \{0\},\quad
  \boldsymbol v = \boldsymbol{0}, \text{ on } \partial \Omega\times (0, T]\,.
\]
The initial space-time mesh $\mathcal{T}_{h}\otimes\mathcal{M}_{\tau}$ is a uniform
tensor-product partition of the space-time domain $\Omega\times I$. We use discretizations
with varying polynomial degrees $r \in \{0,\dots,\,3\}$ in time, where we use
\(\mathbb Q_{r+1}/\mathbb P_{r}^{\mathrm{disc}}\) elements in space, except for
the $dG(0)$ case where we use \(\mathbb Q_{2}/\mathbb P_{1}^{\mathrm{disc}}\) elements.
The goal quantity is chosen to control the $L^2$-error
$e_N^-, e_N^-=u(\boldsymbol{x},T)-u_{\tau h}(T^{-})$, at the final time point $T=1$, given by
\begin{equation}
\label{eq:JL2final}
J_T(u)= \frac{\big(u(\boldsymbol{x},T),e_N^-\big)}{\|e_N^-\|}\,.
\end{equation}
In Tables~\ref{tab:dG0}--\ref{tab:dG3} we present the development of the total discretization error $J(e)=\|e_T\|$, the approximated spatial and temporal error estimators $\eta_h$ and $\eta_\tau$ as well as the effectivity index $\Ieff$ during an adaptive refinement process for various discretizations. For all discretizations considered, we observe a consistent and satisfactory reduction of the discretization error $\|e_T\|$.
Concerning the accuracy of the employed error estimators, all cases exhibit reliable quantitative agreement with the discretization error, as indicated by effectivity indices $\Ieff$ close to one (cf. the last column of Tables~\ref{tab:dG0}--\ref{tab:dG3}).
From the perspective of computational efficiency in space–time adaptive algorithms, it is crucial to ensure an equilibrated reduction of the temporal as well as spatial discretization error.
Referring to this, we point out for all test cases a good equilibration of the spatial and temporal error indicators $\eta_h$ and $\eta_\tau$ in the course of the refinement process (cf. column 6 and 7 of Tables~\ref{tab:dG0}--\ref{tab:dG3}).
\begin{table}[H]\setlength{\tabcolsep}{5.5pt}
\caption{\label{tab:dG0}
Adaptive refinement results for the dG($0$)-$Q_2/P_1^{\mathrm{disc}}$ discretization using the goal quantity \eqref{eq:JL2final} for Example~1.}
\begin{center}\scriptsize
\setlength{\tabcolsep}{8pt}
\resizebox{.99\linewidth}{!}{%
\begin{tabular}
{S[scientific-notation=false,round-precision=0,table-format=1]
    S[scientific-notation=false,round-precision=0,table-format=2]
    S[scientific-notation=false,round-precision=0,table-format=6]
    S[scientific-notation=false,round-precision=0,table-format=7]
    S[table-format=1.1e2]
    S[table-format=1.1e2]
    S[table-format=1.1e2]
    S[table-format=1.1e2]
    S[scientific-notation=false]}
\toprule
\mc{$\ell$} & \mc{$N_{\text{time}}$}
& \mc{$N_{\text{space}}^{\text{max}}$}
& \mc{$N_{\text{tot}}$}
& \mc{$\|e_T\|$}
& \mc{$\tilde{\eta}_{h}$} & \mc{$\tilde{\eta}_\tau$}
& \mc{$\tilde{\eta}_{\tau h}$} & \mc{$\Ieff$}
\\
  \midrule
1        &5      &210       &1050      &2.2017e-02  &1.9422e-03 &1.0825e-02 &1.2767e-02 &0.57989\\
2        &6      &374       &2680      &1.0655e-02  &4.6144e-03 &6.2266e-03 &1.0841e-02 &1.0175\\
3        &7      &592       &4766      &5.1732e-03  &1.1152e-03 &4.0802e-03 &5.1954e-03 &1.0043\\
4        &9      &623       &11069     &2.6679e-03  &1.5040e-03 &2.1946e-03 &3.6986e-03 &1.3863\\
5        &11     &703       &21704     &1.5606e-03  &1.0923e-03 &1.4662e-03 &2.5585e-03 &1.6395\\
6        &14     &761       &40028     &9.7239e-04  &6.4904e-04 &9.5853e-04 &1.6076e-03 &1.6532\\
7        &18     &11761     &91930     &6.4359e-04  &3.4885e-04 &6.6339e-04 &1.0122e-03 &1.5728\\
8        &23     &15203     &193063    &4.4344e-04  &1.4049e-04 &4.4870e-04 &5.8919e-04 &1.3287\\
9        &29     &15489     &439491    &3.2443e-04  &8.0905e-05 &3.3948e-04 &4.2039e-04 &1.2958\\
10       &37     &15187     &869234    &2.3714e-04  &6.5660e-05 &2.4683e-04 &3.1249e-04 &1.3177\\
11       &48     &17776     &1698684   &1.7576e-04  &4.3620e-05 &1.8111e-04 &2.2473e-04 &1.2786\\
12       &62     &26742     &3599510   &1.2984e-04  &3.5041e-05 &1.3323e-04 &1.6827e-04 &1.2960\\
13       &80     &21531     &7238636   &9.8363e-05  &2.2881e-05 &1.0141e-04 &1.2430e-04 &1.2636\\
\bottomrule
\end{tabular}
}
\end{center}
\end{table}
\begin{table}[H]\setlength{\tabcolsep}{4.5pt}
\caption{\label{tab:dG1}
Adaptive refinement results for the dG($1$)-$Q_2/P_1^{\mathrm{disc}}$ discretization using the goal quantity \eqref{eq:JL2final} for Example~1.
}
\begin{center}\scriptsize
\setlength{\tabcolsep}{8pt}
\resizebox{.99\linewidth}{!}{%
\begin{tabular}
{S[scientific-notation=false,round-precision=0,table-format=1]
    S[scientific-notation=false,round-precision=0,table-format=2]
    S[scientific-notation=false,round-precision=0,table-format=6]
    S[scientific-notation=false,round-precision=0,table-format=7]
    S[table-format=1.1e2]
    S[table-format=1.1e2]
    S[table-format=1.1e2]
    S[table-format=1.1e2]
    S[scientific-notation=false]}
\toprule
\mc{$\ell$} & \mc{$N_{\text{time}}$}
& \mc{$N_{\text{space}}^{\text{max}}$}
& \mc{$N_{\text{tot}}$}
& \mc{$\|e_T\|$}
& \mc{$\tilde{\eta}_{h}$} & \mc{$\tilde{\eta}_\tau$}
& \mc{$\tilde{\eta}_{\tau h}$} & \mc{$\Ieff$}
\\
\midrule
1  &   4  &   1540   &  6160    &9.3412e-04 &2.8871e-04  &6.7677e-05 &3.5639e-04  &0.38\\
2  &   5  &   1736   &  12390   &4.5738e-04 &3.1714e-04  &4.0001e-05 &3.5714e-04  &0.78\\
3  &   5  &   3248   &  23138   &2.1601e-04 &-7.3945e-05 &1.9074e-05 &-5.4871e-05 &0.25\\
4  &   7  &   3494   &  57272   &9.1844e-05 &-4.3569e-06 &4.4129e-06 &5.5998e-08  &0.00061\\
5  &   9  &   3572   &  114154  &5.1669e-05 &5.9771e-05  &6.6503e-06 &6.6421e-05  &1.29\\
6  &   9  &   4712   &  183120  &2.5171e-05 &8.6274e-06  &2.3658e-06 &1.0993e-05  &0.44\\
7  &  12  &   4868   &  371674  &1.0561e-05 &9.7207e-06  &3.2539e-07 &1.0046e-05  &0.95\\
8  &  12  &   4868   &  655942  &6.2456e-06 &3.8154e-06  &4.4631e-07 &4.2617e-06  &0.68\\
9  &  12  &   129014 &  1180326 &2.7837e-06 &-3.4125e-06 &3.9546e-07 &-3.0170e-06 &1.08\\
\bottomrule
\end{tabular}
}
\end{center}
\end{table}

\begin{figure}[H]
\centering
\begin{tikzpicture}
\begin{groupplot}[
group style={
  group size=3 by 1,
  horizontal sep=1.5cm,
},
width=0.34\textwidth,
height=0.32\textwidth,
grid=both,
xlabel={adaptive loop $\ell$},
legend columns=4]

\nextgroupplot[
ylabel={$N_{\mathrm{Newt}}/N_{\mathrm{slabs}}$},
legend style={/tikz/every even column/.append style={column sep=0.3cm}, draw=none, fill=none, font=\footnotesize, at={(-0.2,1.1)}, anchor=south west}]
\addplot+[color=black,mark=*,mark options={fill=black}] table[x=loop,y=Newton_per_timestep,col sep=comma] {DG0_work_table.csv};
\addlegendentry{$dG(0)$-$Q_2/P_1^{\mathrm{disc}}$}
\addplot+[color=orange,mark=square*,mark options={fill=orange}] table[x=loop,y=Newton_per_timestep,col sep=comma] {DG1_work_table.csv};
\addlegendentry{$dG(1)$-$Q_2/P_1^{\mathrm{disc}}$}
\addplot+[color=HSUred,mark=triangle*,mark options={fill=HSUred}] table[x=loop,y=Newton_per_timestep,col sep=comma] {DG2_work_table.csv};
\addlegendentry{$dG(2)$-$Q_3/P_2^{\mathrm{disc}}$}
\addplot+[color=blue,mark=diamond*,mark options={fill=blue}] table[x=loop,y=Newton_per_timestep,col sep=comma] {DG3_work_table.csv};
\addlegendentry{$dG(3)$-$Q_4/P_3^{\mathrm{disc}}$}

\nextgroupplot[
ylabel={$N_{\mathrm{lin}}^{p}/N_{\mathrm{Newt}}$},
]
\addplot+[color=black,mark=*,mark options={fill=black}] table[x=loop,y=LinearPrimal_per_Newton,col sep=comma] {DG0_work_table.csv};
\addplot+[color=orange,mark=square*,mark options={fill=orange}] table[x=loop,y=LinearPrimal_per_Newton,col sep=comma] {DG1_work_table.csv};
\addplot+[color=HSUred,mark=triangle*,mark options={fill=HSUred}] table[x=loop,y=LinearPrimal_per_Newton,col sep=comma] {DG2_work_table.csv};
\addplot+[color=blue,mark=diamond*,mark options={fill=blue}] table[x=loop,y=LinearPrimal_per_Newton,col sep=comma] {DG3_work_table.csv};

\nextgroupplot[
ylabel={$W_d/W_{\mathrm{tot}}$},
ymin=0,
ymax=0.5,
]
\addplot+[color=black,mark=*,mark options={fill=black}] table[x=loop,y expr=\thisrow{Work_dual}/\thisrow{Work_total},col sep=comma] {DG0_work_table.csv};
\addplot+[color=orange,mark=square*,mark options={fill=orange}] table[x=loop,y expr=\thisrow{Work_dual}/\thisrow{Work_total},col sep=comma] {DG1_work_table.csv};
\addplot+[color=HSUred,mark=triangle*,mark options={fill=HSUred}] table[x=loop,y expr=\thisrow{Work_dual}/\thisrow{Work_total},col sep=comma] {DG2_work_table.csv};
\addplot+[color=blue,mark=diamond*,mark options={fill=blue}] table[x=loop,y expr=\thisrow{Work_dual}/\thisrow{Work_total},col sep=comma] {DG3_work_table.csv};

\end{groupplot}
\end{tikzpicture}
\caption{Algebraic performance over the adaptive loops: nonlinear iterations per slab, primal linear iterations per nonlinear step, and relative dual work fraction.}
\label{fig:example1-solver-groupplot}
\end{figure}
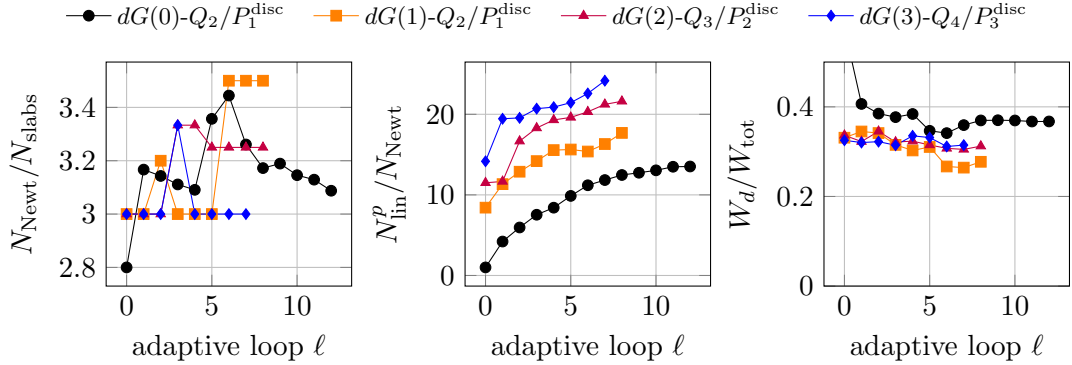

\begin{figure}[H]
\centering
  \begin{tikzpicture}
    \pgfplotsset{every axis/.style={ ymin=1000, ymax=1000000,
        xmin=-0.5,xmax=2.5,%
        width=0.5\linewidth,height=0.4\linewidth,name=rt3d,%
        xtick={0,0.66,1.33,2}, xticklabels={{\scriptsize$Q_2/P_1^{\mathrm{disc}}$},
        {\scriptsize$Q_2/P_1^{\mathrm{disc}}$}, {\scriptsize$Q_3/P_2^{\mathrm{disc}}$}, {\scriptsize$Q_4/P_3^{\mathrm{disc}}$}},%
        xlabel=FE pair,ybar, ymajorgrids,%
        major grid style={white,thick}, clip=false,%
        axis on top, log origin=infty,%
        legend style={%
          /tikz/every even column/.append style={column sep=0.7cm},%
          at={(0.35,1.1)},%
          anchor=south,%
          font=\footnotesize, draw=none},%
        bar width=16pt}%
    }%
    \begin{groupplot}[group style={columns=2,group name=rte,horizontal
        sep=.12\textwidth},%
      width=0.5\textwidth,%
      height=0.36\textwidth, axis line style=thick]
      \nextgroupplot[bar shift=0pt,stack plots=y,ymode=log,ymin=1000,
      ymax=1000000, legend columns=-1,ylabel={time $[\text{s}]$},%
      axis on top]%
      \addlegendimage{draw=blue,fill=blue}%
      \addlegendimage{draw=red,fill=red}
      \addlegendimage{draw=green,fill=green}
      \addlegendentry{primal}%
      \addlegendentry{dual}%
      \addlegendentry{error estimator}%
      \addplot[draw=none,fill=blue,forget plot] coordinates {%
        (0, 3550) (0.66, 2600) (1.33, 57100) (2, 80900)%
      };%
      \addplot[draw=none,fill=red,forget plot] coordinates {%
        (0, 5510) (0.66, 3780) (1.33, 85800) (2, 121100)%
      };%
      \addplot[draw=none,fill=green,forget plot] coordinates {%
        (0, 5612) (0.66, 3808) (1.33, 86099) (2, 121375)%
      };%
      \nextgroupplot[bar shift=0pt,stack plots=y,ymin=0,legend
      columns=3,ylabel={Portion of wall time $[\%]$}, ymax=100]%
      \addlegendimage{draw=blue,fill=blue}%
      \addlegendimage{draw=red,fill=red}
      \addlegendimage{draw=green,fill=green}
      \addlegendentry{primal}%
      \addlegendentry{dual}%
      \addlegendentry{error estimator}%
      \addplot[draw=none,fill=blue,forget plot] coordinates {%
        (0, 64.00) (0.66, 68.0) (1.33, 66.0) (2, 66.0)%
      };%
      \addplot[draw=none,fill=red,forget plot] coordinates {%
        (0, 34.00) (0.66, 31.0) (1.33, 33.0) (2, 33.0)%
      };%
      \addplot[draw=none,fill=green,forget plot] coordinates {%
        (0, 2.00) (0.66, 1.0) (1.33, 1.0) (2, 1.0)%
      };%
    \end{groupplot}
    \node at (rte c1r1.north west)[anchor=north west] {(a)};%
    \node at (rte c2r1.north west)[anchor=north west] {(b)};%
  \end{tikzpicture}
\caption{Share of time spent on the primal problem, dual problem, and error estimation over the total wall-clock time for an adaptive run, including all DWR loops, measured in seconds (left), and their percentage shares (right) for the different discretizations corresponding to Tables~\ref{tab:dG0}--\ref{tab:dG3}.
}
\label{fig:walltimes}
\end{figure}
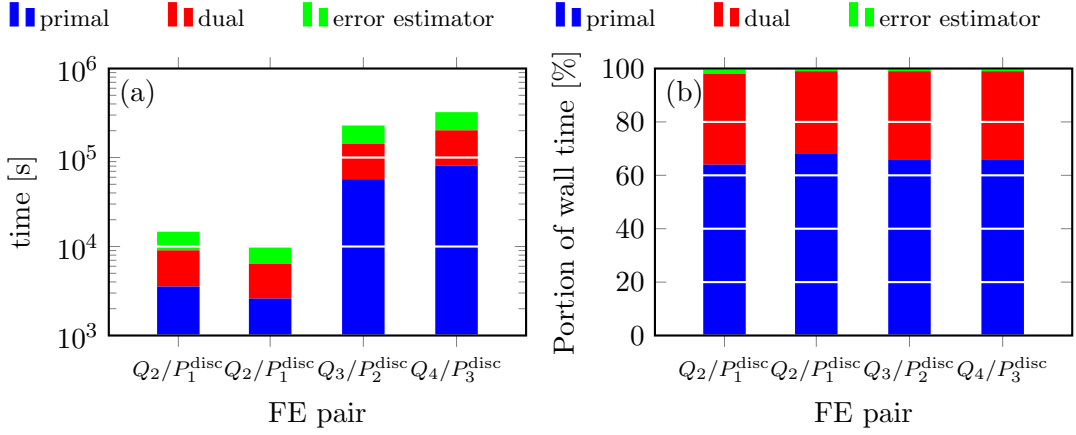

\begin{table}[H]\setlength{\tabcolsep}{4.5pt}
\caption{\label{tab:dG2}
Adaptive refinement results for the dG($2$)-$Q_3/P_2^{\mathrm{disc}}$ discretization using the goal quantity \eqref{eq:JL2final} for Example~1.
}
\begin{center}\scriptsize
\setlength{\tabcolsep}{8pt}
\resizebox{.99\linewidth}{!}{%
\begin{tabular}
{S[scientific-notation=false,round-precision=0,table-format=1]
    S[scientific-notation=false,round-precision=0,table-format=2]
    S[scientific-notation=false,round-precision=0,table-format=6]
    S[scientific-notation=false,round-precision=0,table-format=7]
    S[table-format=1.1e2]
    S[table-format=1.1e2]
    S[table-format=1.1e2]
    S[table-format=1.1e2]
    S[scientific-notation=false]}
\toprule
\mc{$\ell$} & \mc{$N_{\text{time}}$}
& \mc{$N_{\text{space}}^{\text{max}}$}
& \mc{$N_{\text{tot}}$}
& \mc{$\|e_{\text{drag}}\|$}
& \mc{$\tilde{\eta}_{h}$} & \mc{$\tilde{\eta}_\tau$}
& \mc{$\tilde{\eta}_{\tau h}$} & \mc{$\Ieff$}
\\
\midrule
1  &  2  &    4902  &     9804 & 1.7087e-04 &  2.9316e-06 &4.3518e-05 &4.6450e-05 &0.27184\\
2  &  3  &    4902  &    14706 & 3.4100e-05 &  1.5325e-05 &1.5867e-06 &1.6912e-05 &0.49595\\
3  &  3  &   19920  &    61122 & 1.9126e-05 &  6.5756e-05 &2.9991e-06 &6.8755e-05 &3.5949\\
4  &  3  &   40368  &   128460 & 1.6319e-05 &  4.7500e-05 &3.5835e-06 &5.1084e-05 &3.1303\\
5  &  3  &  121374  &   399990 & 1.4557e-05 &  1.4939e-05 &4.0396e-06 &1.8979e-05 &1.3038\\
6  &  4  &  137706  &   782160 & 1.0120e-06 &  6.7099e-06 &1.8632e-07 &6.8962e-06 &6.8141\\
7  &  4  &  234618  &  1229250 & 9.5143e-07 & 3.5681e-06  &1.9791e-07 &3.7660e-06 &3.9583\\
8 &  4  &  704400  &  3852132 & 9.1104e-07 & 1.3522e-06  &2.0716e-07 &1.5593e-06 &1.7116\\
9 &  4  & 1106778  &  6574092 & 9.0892e-07 & 8.4627e-07  &2.0766e-07 &1.0539e-06 &1.1595\\
\bottomrule
\end{tabular}
}
\end{center}
\end{table}
\begin{table}[H]\setlength{\tabcolsep}{4.5pt}
\caption{\label{tab:dG3}
Adaptive refinement results for the dG($3$)-$Q_4/P_3^{\mathrm{disc}}$ discretization using the goal quantity \eqref{eq:JL2final} for Example~1.
}
\begin{center}\scriptsize
\setlength{\tabcolsep}{8pt}
\resizebox{.99\linewidth}{!}{%
\begin{tabular}
{S[scientific-notation=false,round-precision=0,table-format=1]
    S[scientific-notation=false,round-precision=0,table-format=2]
    S[scientific-notation=false,round-precision=0,table-format=6]
    S[scientific-notation=false,round-precision=0,table-format=7]
    S[table-format=1.1e2]
    S[table-format=1.1e2]
    S[table-format=1.1e2]
    S[table-format=1.1e2]
    S[scientific-notation=false]}
\toprule
\mc{$\ell$} & \mc{$N_{\text{time}}$}
& \mc{$N_{\text{space}}^{\text{max}}$}
& \mc{$N_{\text{tot}}$}
& \mc{$\|e_T\|$}
& \mc{$\tilde{\eta}_{h}$} & \mc{$\tilde{\eta}_\tau$}
& \mc{$\tilde{\eta}_{\tau h}$} & \mc{$\Ieff$}
\\
\midrule
1   &2 &11272     &22544   &4.6212e-06  &1.8080e-07  &4.3643e-07 &6.1723e-07  &0.13\\
2   &3 &18664     &78888   &7.6425e-07  &-5.3932e-07 &5.3089e-09 &-5.3401e-07 &0.70\\
3   &3 &38536     &157424  &3.3388e-07  &-4.0245e-07 &1.4767e-08 &-3.8769e-07 &1.16\\
4   &3 &53528     &277208  &1.7752e-07  &-2.6782e-07 &2.2228e-08 &-2.4559e-07 &1.38\\
5   &3 &62944     &450472  &1.7548e-07  &-2.9604e-07 &2.2457e-08 &-2.7359e-07 &1.56\\
6   &3 &113528    &596728  &1.7485e-07  &-1.9283e-07 &2.2517e-08 &-1.7032e-07 &0.97\\
7   &4 &225000    &1958424 &6.1465e-09 &-6.62e-09 &4.36e-10 &-6.1870e-09 &1.01\\
8   &4 &284728    &2768280 &4.3765e-09 &-6.89e-09 &6.13e-10 &-6.2785e-09 &1.43\\
\bottomrule
\end{tabular}
}
\end{center}
\end{table}
Across all four fully adaptive runs, the DWR estimator in combination with
adaptive refinement yields a consistent reduction of the goal error with
increasing space-time resolution, while the corresponding effectivity indices
remain in a reasonable range. The reported mean rates show that the decay of the
goal error with respect to total work is close to its decay with respect to the
space-time number of degrees of freedom, indicating that the adaptive refinement
translates effectively into computational efficiency. A compact comparison of
these approximation and cost quantities is given in Table~\ref{tab:dwr-summary}.
\begin{table}[H]
  \centering
  \caption{\label{tab:dwr-summary}Summary of the fully adaptive DWR computations for the four temporal
    orders \(r\). Reported are the data from the final adaptive cycle: The
    space-time number of degrees of freedom \(N_{\mathrm{ST}}\), the goal error
    \(\|e_T\|\), the magnitude of \(|\tilde{\eta}_{\tau h}|\), and
    the effectivity index. The table includes the convergence
    rates~\eqref{eq:mean-rates}, as well as the ranges of the average nonlinear
    iterations, the average primal linear iterations per nonlinear step~\eqref{eq:iteration-indicators}, and the relative dual work fraction~\eqref{eq:dual-work-fraction}.}
\footnotesize
\resizebox{.99\linewidth}{!}{%
\begin{tabular}{
  S[scientific-notation=false,round-precision=0,table-format=1]
  S[scientific-notation=false,round-precision=0,table-format=2]
  S[scientific-notation=false,round-precision=0,table-format=6]
  |
  S[table-format=1.2e2]
  S[table-format=1.2e2]
  S[scientific-notation=false]
  |
  S[scientific-notation=false]
  S[scientific-notation=false]
  |
  S[scientific-notation=false]
  S[scientific-notation=false]
  |
  S[scientific-notation=false]
  S[scientific-notation=false]
  |
  S[scientific-notation=false]
  }
  \toprule
  &
  &
  &\multicolumn{3}{c|}{DWR}
  &\multicolumn{2}{c|}{mean rates}
  &\multicolumn{2}{c|}{$n_{\mathrm{NL}}/N_{\mathrm{time}}$}
  &\multicolumn{2}{c|}{$n_{\mathrm{L}}^{p}/n_{\mathrm{NL}}$}
  \\
  \mc{$r$}
  & \mc{$\ell_{\max}$}
  & \multicolumn{1}{c|}{$N_{\mathrm{ST}}$}
  & \mc{$\| e_T \|$}
  & \mc{$|\tilde{\eta}_{\tau h}|$}
  & \multicolumn{1}{c|}{$I_{\mathrm{eff}}$}
  & \mc{$\overline{\operatorname{rate}}_{N_{\mathrm{ST}}}^{\|e_T\|}$}
  & \multicolumn{1}{c|}{$\overline{\operatorname{rate}}_{W_{\mathrm{tot}}}^{\|e_T\|}$}
  & \mc{$\min$}
  & \multicolumn{1}{c|}{$\max$}
  & \mc{$\min$}
  & \multicolumn{1}{c|}{$\max$}
  & \mc{$W_d/W_{\mathrm{tot}}$} \\
\midrule
0 & 13 & 7238636 & 9.84e-05 & 1.24e-04 & 1.26 & 0.62 & 0.50 & 2.8 & 3.44 & 1.0 & 13.52 & 0.37
\\
1 & 9  & 1180326 & 2.78e-06 & 3.02e-06 & 1.08 & 1.13 & 0.99 & 3.0 & 3.50 & 8.42 & 17.69 & 0.28
\\
2 & 9  & 6574092 & 9.09e-07 & 1.05e-06 & 1.16 & 1.11 & 1.12 & 3.0 & 3.33 & 11.5 & 21.62 & 0.31
\\
3 & 8  & 2768280 & 4.38e-09 & 6.28e-09 & 1.43 & 1.08 & 0.98 & 3.0 & 3.33 & 14.17 & 24.17 & 0.31
\\
\bottomrule
\end{tabular}
}
\end{table}
The algebraic behavior is summarized in
Figure~\ref{fig:example1-solver-groupplot}, which shows stable nonlinear
iteration counts per slab and only moderate growth of the primal linear
iterations per nonlinear step over the adaptive loops. This indicates that the
slab-wise geometric multigrid preconditioned Newton--Krylov method remains
robust under successive adaptive refinement. At the same time, the relative dual
work fraction stays at a comparable level across all discretizations, showing
that the adjoint solve contributes a noticeable but controlled part of the overall computational cost.
In Fig.~\ref{fig:walltimes} we plot the share of the primal problem, the dual problem as well as the error estimation over the entire wallclock time and its percentage share for one adaptive run including all DWR loops for the different discretizations corresponding to Tables~\ref{tab:dG0}--\ref{tab:dG3}.
While the total runtime from dG($1$)-$Q_2/P_1^{\textnormal{disc}}$ to dG($3$)-$Q_4/P_3^{\textnormal{disc}}$ increases progressively as expected, the percentage share of the three units remains stable and nearly unchanged across all discretizations. In comparison to the primal problem, the dual problem requires only one third of its computational effort, while the error estimator itself is negligible at less than one percent.
The fact that the total runtime of the dG($0$)-$Q_2/P_1^{\textnormal{disc}}$ discretization is slightly higher to that of the dG($1$)-$Q_2/P_1^{\textnormal{disc}}$ one is due to the large number of slabs required in the first case.

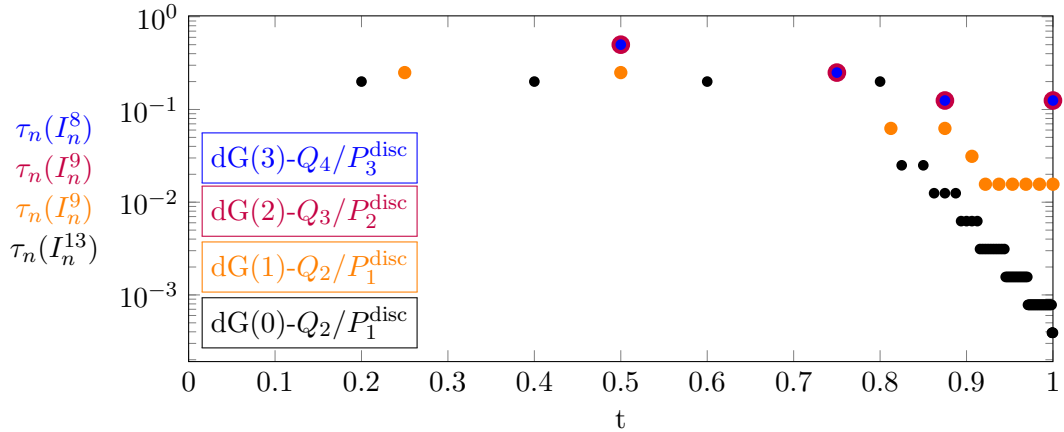
\begin{figure}[htb]
\begin{minipage}{\linewidth}
\centering
\begin{tikzpicture}
\begin{axis}[%
width=4.5in,
height=1.8in,
scale only axis,
/pgf/number format/.cd, 1000 sep={},
xlabel={t},
ylabel style={rotate=-90},
ylabel style={align=center},
ylabel={
\textcolor{blue}{$\tau_n(I_n^{8})$}\\[1.5pt]
\textcolor{HSUred}{$\tau_n(I_n^{9})$}\\[1.5pt]
\textcolor{orange}{$\tau_n(I_n^{9})$}\\[1.5pt]
$\tau_n(I_n^{13})$
},
ymode=log,
xmin=0.0,
xmax=1.0,
yminorticks=true,
]
\node at (axis cs:0.14,5e-4) {\fbox{dG($0$)-$Q_2/P_1^{\textnormal{disc}}$}};
\node at (axis cs:0.14,2e-3) {\textcolor{orange}{\fbox{dG($1$)-$Q_2/P_1^{\textnormal{disc}}$}}};
\node at (axis cs:0.14,8e-3) {\textcolor{HSUred}{\fbox{dG($2$)-$Q_3/P_2^{\textnormal{disc}}$}}};
\node at (axis cs:0.14,3e-2) {\textcolor{blue}{\fbox{dG($3$)-$Q_4/P_3^{\textnormal{disc}}$}}};
\addplot [
color=black,
solid,
line width=1.0pt,
mark=*,
mark size = 1.5,
only marks,
mark options={solid,black}
]
table[row sep=crcr]{
0.2 0.2 \\
0.4 0.2 \\
0.6 0.2 \\
0.8 0.2 \\
0.825 0.025 \\
0.85 0.025 \\
0.8625 0.0125 \\
0.875 0.0125 \\
0.8875 0.0125 \\
0.89375 0.00625 \\
0.9 0.00625 \\
0.90625 0.00625 \\
0.9125 0.00625 \\
0.915625 0.003125 \\
0.91875 0.003125 \\
0.921875 0.003125 \\
0.925 0.003125 \\
0.928125 0.003125 \\
0.93125 0.003125 \\
0.934375 0.003125 \\
0.9375 0.003125 \\
0.940625 0.003125 \\
0.94375 0.003125 \\
0.945312 0.0015625 \\
0.946875 0.0015625 \\
0.948437 0.0015625 \\
0.95 0.0015625 \\
0.951562 0.0015625 \\
0.953125 0.0015625 \\
0.954687 0.0015625 \\
0.95625 0.0015625 \\
0.957812 0.0015625 \\
0.959375 0.0015625 \\
0.960937 0.0015625 \\
0.9625 0.0015625 \\
0.964062 0.0015625 \\
0.965625 0.0015625 \\
0.967187 0.0015625 \\
0.96875 0.0015625 \\
0.970313 0.0015625 \\
0.971094 0.00078125 \\
0.971875 0.00078125 \\
0.972656 0.00078125 \\
0.973437 0.00078125 \\
0.974219 0.00078125 \\
0.975 0.00078125 \\
0.975781 0.00078125 \\
0.976562 0.00078125 \\
0.977344 0.00078125 \\
0.978125 0.00078125 \\
0.978906 0.00078125 \\
0.979687 0.00078125 \\
0.980469 0.00078125 \\
0.98125 0.00078125 \\
0.982031 0.00078125 \\
0.982812 0.00078125 \\
0.983594 0.00078125 \\
0.984375 0.00078125 \\
0.985156 0.00078125 \\
0.985938 0.00078125 \\
0.986719 0.00078125 \\
0.9875 0.00078125 \\
0.988281 0.00078125 \\
0.989063 0.00078125 \\
0.989844 0.00078125 \\
0.990625 0.00078125 \\
0.991406 0.00078125 \\
0.992188 0.00078125 \\
0.992969 0.00078125 \\
0.99375 0.00078125 \\
0.994531 0.00078125 \\
0.995313 0.00078125 \\
0.996094 0.00078125 \\
0.996875 0.00078125 \\
0.997656 0.00078125 \\
0.998437 0.00078125 \\
0.998828 0.000390625 \\
0.999219 0.000390625 \\
0.999609 0.000390625 \\
1 0.000390625 \\
};
\addplot [
color=black,
solid,
line width=1.0pt,
mark=*,
mark size = 2.,
only marks,
mark options={solid,orange}
]
table[row sep=crcr]{
0.25 0.25 \\
0.5 0.25 \\
0.75 0.25 \\
0.8125 0.0625 \\
0.875 0.0625 \\
0.90625 0.03125 \\
0.921875 0.015625 \\
0.9375 0.015625 \\
0.953125 0.015625 \\
0.96875 0.015625 \\
0.984375 0.015625 \\
1 0.015625 \\
};
\addplot [
color=black,
solid,
line width=1.0pt,
mark=*,
mark size = 3.0,
only marks,
mark options={solid,HSUred}
]
table[row sep=crcr]{
0.5 0.5 \\
0.75 0.25\\
0.875 0.125\\
1 0.125 \\
};
\addplot [
color=black,
solid,
line width=1.0pt,
mark=*,
mark size = 1.5,
only marks,
mark options={solid,blue}
]
table[row sep=crcr]{
0.5 0.5 \\
0.75 0.25\\
0.875 0.125\\
1 0.125 \\
};
\end{axis}
\end{tikzpicture}
\end{minipage}
\caption{Distribution of the adaptively determined temporal step size $\tau_n$ plotted over the entire time interval $I=(0,T]$ for the different discretizations corresponding to the respective final DWR loop in
Tables~\ref{tab:dG0}--\ref{tab:dG3}.
}
\label{BruchhaeuserFig:tau-dist}
\end{figure}

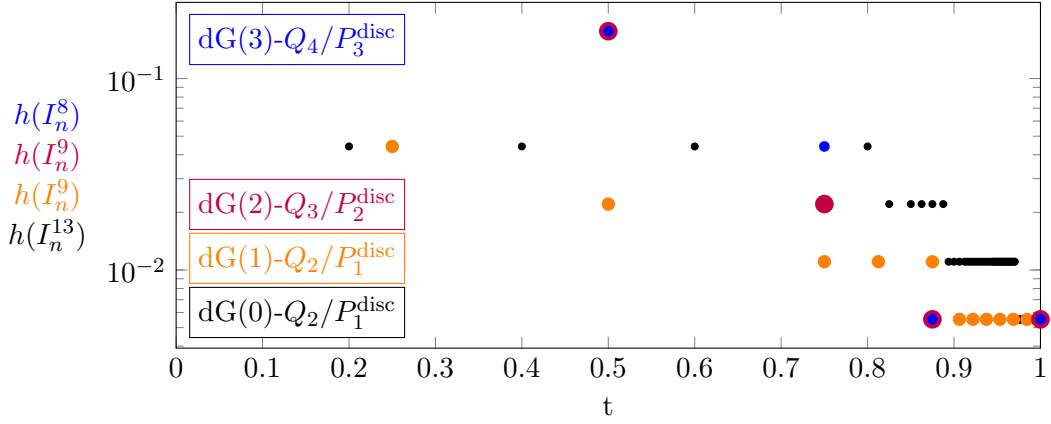
\begin{figure}[htb]
\begin{minipage}{\linewidth}
\centering
\begin{tikzpicture}
\begin{axis}[%
width=4.5in,
height=1.8in,
scale only axis,
/pgf/number format/.cd, 1000 sep={},
xlabel={t},
ylabel style={rotate=-90},
ylabel style={align=center},
ylabel={
\textcolor{blue}{$h(I_n^{8})$}\\[1.5pt]
\textcolor{HSUred}{$h(I_n^{9})$}\\[1.5pt]
\textcolor{orange}{$h(I_n^{9})$}\\[1.5pt]
$h(I_n^{13})$
},
xmin=0.0,
xmax=1.0,
ymode=log,
yminorticks=true,
]
\node at (axis cs:0.14,6e-3) {\fbox{dG($0$)-$Q_2/P_1^{\textnormal{disc}}$}};
\node at (axis cs:0.14,1.15e-2) {\textcolor{orange}{\fbox{dG($1$)-$Q_2/P_1^{\textnormal{disc}}$}}};
\node at (axis cs:0.14,2.2e-2) {\textcolor{HSUred}{\fbox{dG($2$)-$Q_3/P_2^{\textnormal{disc}}$}}};
\node at (axis cs:0.14,1.6e-1) {\textcolor{blue}
{\fbox{dG($3$)-$Q_4/P_3^{\textnormal{disc}}$}}};
\addplot [
color=black,
solid,
line width=1.0pt,
mark=*,
mark size = 1.,
only marks,
mark options={solid,black}
]
table[row sep=crcr]{
0.2 0.0441942 \\
0.4 0.0441942 \\
0.6 0.0441942 \\
0.8 0.0441942 \\
0.825 0.0220971 \\
0.85 0.0220971 \\
0.8625 0.0220971 \\
0.875 0.0220971 \\
0.8875 0.0220971 \\
0.89375 0.0110485 \\
0.9 0.0110485 \\
0.90625 0.0110485 \\
0.9125 0.0110485 \\
0.915625 0.0110485 \\
0.91875 0.0110485 \\
0.921875 0.0110485 \\
0.925 0.0110485 \\
0.928125 0.0110485 \\
0.93125 0.0110485 \\
0.934375 0.0110485 \\
0.9375 0.0110485 \\
0.940625 0.0110485 \\
0.94375 0.0110485 \\
0.945312 0.0110485 \\
0.946875 0.0110485 \\
0.948437 0.0110485 \\
0.95 0.0110485 \\
0.951562 0.0110485 \\
0.953125 0.0110485 \\
0.954687 0.0110485 \\
0.95625 0.0110485 \\
0.957812 0.0110485 \\
0.959375 0.0110485 \\
0.960937 0.0110485 \\
0.9625 0.0110485 \\
0.964062 0.0110485 \\
0.965625 0.0110485 \\
0.967187 0.0110485 \\
0.96875 0.0110485 \\
0.970313 0.0110485 \\
0.971094 0.00552427 \\
0.971875 0.00552427 \\
0.972656 0.00552427 \\
0.973437 0.00552427 \\
0.974219 0.00552427 \\
0.975 0.00552427 \\
0.975781 0.00552427 \\
0.976562 0.00552427 \\
0.977344 0.00552427 \\
0.978125 0.00552427 \\
0.978906 0.00552427 \\
0.979687 0.00552427 \\
0.980469 0.00552427 \\
0.98125 0.00552427 \\
0.982031 0.00552427 \\
0.982812 0.00552427 \\
0.983594 0.00552427 \\
0.984375 0.00552427 \\
0.985156 0.00552427 \\
0.985938 0.00552427 \\
0.986719 0.00552427 \\
0.9875 0.00552427 \\
0.988281 0.00552427 \\
0.989063 0.00552427 \\
0.989844 0.00552427 \\
0.990625 0.00552427 \\
0.991406 0.00552427 \\
0.992188 0.00552427 \\
0.992969 0.00552427 \\
0.99375 0.00552427 \\
0.994531 0.00552427 \\
0.995313 0.00552427 \\
0.996094 0.00552427 \\
0.996875 0.00552427 \\
0.997656 0.00552427 \\
0.998437 0.00552427 \\
0.998828 0.00552427 \\
0.999219 0.00552427 \\
0.999609 0.00552427 \\
1 0.00552427 \\
};
\addplot [
color=black,
solid,
line width=1.0pt,
mark=*,
mark size = 2.,
only marks,
mark options={solid,orange}
]
table[row sep=crcr]{
0.25 0.0441942 \\
0.5 0.0220971 \\
0.75 0.0110485 \\
0.8125 0.0110485 \\
0.875 0.0110485 \\
0.90625 0.00552427 \\
0.921875 0.00552427 \\
0.9375 0.00552427 \\
0.953125 0.00552427 \\
0.96875 0.00552427 \\
0.984375 0.00552427 \\
1 0.00552427 \\
};
\addplot [
color=black,
solid,
line width=1.0pt,
mark=*,
mark size = 3.0,
only marks,
mark options={solid,HSUred}
]
table[row sep=crcr]{
0.5 0.176777 \\
0.75 0.0220971 \\
0.875 0.00552427 \\
1 0.00552427 \\
};
\addplot [
color=black,
solid,
line width=1.0pt,
mark=*,
mark size = 1.5,
only marks,
mark options={solid,blue}
]
table[row sep=crcr]{
0.5 0.176777 \\
0.75 0.0441942 \\
0.875 0.00552427 \\
1 0.00552427 \\
};
\end{axis}
\end{tikzpicture}
\end{minipage}
\caption{Distribution of the adaptively determined minimal spatial element diameter $\min_{K\in\mathcal{T}_h^{n,\ell_{\max}}} h_K$ within the slab $Q_n$ plotted over the entire time interval $I=(0,T]$ for the different discretizations corresponding to the respective final DWR loop in
Tables~\ref{tab:dG0}--\ref{tab:dG3}.
}
\label{BruchhaeuserFig:h-dist}
\end{figure}

Finally, in Fig.~\ref{BruchhaeuserFig:tau-dist} and Fig.~\ref{BruchhaeuserFig:h-dist} we visualize the distribution of the adaptively determined time step sizes $\tau_n$ and minimal spatial element diameter $\min_{K\in\mathcal{T}_h^{n,\ell_{\max}}} h_K$ within the slab $Q_n$ over the entire time interval $I$, respectively, corresponding to the final DWR loop of the different discretizations listed in Tables~\ref{tab:dG0}--\ref{tab:dG3}.
For all discretizations, the step and element sizes are relatively large at the beginning, but they decrease continuously over time, reaching their minimal values at the final time $T=1$.
We note that the initial steps correspond to the initial distribution of the temporal and spatial meshes.
This behavior is consistent with the choice of the underlying goal functional \eqref{eq:JL2final}, which specifically targets the control of the $L^2$-error at the final time $T$.
Furthermore, we observe that the total length of the subintervals is significantly reduced when using higher-order discretizations compared to the lowest-order case.

\subsection{Example 2: DFG Benchmark 2D-3} %
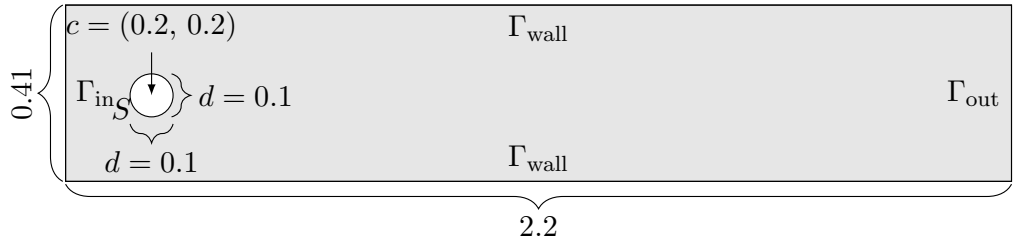
\begin{figure}[H]
  \centering
  \resizebox{.95\linewidth}{!}{%
  \begin{tikzpicture}[scale=5.5]
    \draw [fill=gray!20] (0,0) rectangle (2.2,0.41);
    \path[draw=black] (0,0) -- node[midway, above]{$\Gamma_{\textrm{wall}}$} ++(2.2,0) --
    node[midway, left]{$\Gamma_{\textrm{out}}$} ++(0,0.41) -- node[midway,
    below]{$\Gamma_{\textrm{wall}}$} ++(-2.2,0) -- node[midway, right]{$\Gamma_{\textrm{in}}$}cycle;
    \draw[fill=white] (0.2,0.2) ellipse (0.05 and 0.05);
    \draw[decorate,decoration={brace,raise=1pt,amplitude=9pt,mirror}] (0,0) --
    node[below=9pt]{$2{.}2$} (2.2,0);
    \draw[decorate,decoration={brace,raise=1pt,amplitude=9pt}] (0,0) --
    node[pos=0.5,xshift=-16pt,rotate=90]{$0{.}41$} (0,0.41);
    \draw[decorate,decoration={brace,raise=1pt,amplitude=7pt,mirror}] (0.24,0.15) --
    node[right=7pt]{$d=0{.}1$} (0.24,0.25);
    \draw[decorate,decoration={brace,raise=1pt,amplitude=7pt,mirror}] (0.15,0.14) --
    node[below=7pt]{$d=0{.}1$} (0.25,0.14);
    \draw[-latex] (0.2, 0.3) node[above]{$c=(0{.}2,\,0{.}2)$} --(0.2,0.2);
    \draw (0.13,0.17) node[]{$S$};
  \end{tikzpicture}
  }
\caption{Geometry of the test scenario with a parabolic inflow profile $\Gamma_{\textrm{in}}$,
  do-nothing boundary conditions at the outflow boundary $\Gamma_{\textrm{out}}$ and no-slip conditions on the walls
  $\Gamma_{\textrm{wall}}$.}\label{fig:geom}
\end{figure}

We consider the well-established benchmark problem of laminar flow around a cylinder introduced by Sch\"afer and Turek~\cite{BruchhaeuserST96}, precisely the unsteady 2D-3 test setting. In this case, the flow around a circular obstacle with Reynolds number \(\mathrm{Re}=100\)
(the kinematic viscosity is set to $\nu=0.001$)
is driven by a parabolic inflow profile whose amplitude is modulated in time by a \(\sin\) function. The geometry of the flow domain is shown in Fig.~\ref{fig:geom}. We note that the coarse mesh used here is the same as in previous studies~\cite{BruchhaeuserBR12,RoThiKoeWi23}.
The final time point is set to $T=8$.
The goal quantity is chosen to control the error of the mean drag coeﬃcient. This coefficient is defined by
\begin{equation}
\label{eq:Jmeandrag}
J_{\textnormal{drag}}(u)=
\displaystyle-\int_I\{(\partial_t\velocity,\hat{\primalpsi}_d)+a(\boldsymbol{u})(\hat{\primalphi}_d)\}\,\mathrm{d}t\,,
\end{equation}
where $\hat{\primalphi}_d=(\hat{\primalpsi}_d,0)^\top$ fulfills $\hat{\primalpsi}_d|_S=(|I|^{-1}C,0)^\top\,, \;\; \hat{\primalpsi}_d|_{\partial\Omega\textbackslash S}=\boldsymbol{0}$ with the constant $C=2(U_{\mathrm{mean}}^2d)^{-1}$.
Here, the mean inflow velocity is given by $$U_{\mathrm{mean}}=2/3\;\velocity_D(0,0.205,0)=\sin (\pi t/8)$$ given by Eq.~\eqref{eq:inflow} and $d$ denotes the diameter of the obstacle with boundary $S=S=\partial B_r(0.2,0.2)$.
We take the reference value for the goal functional from~\cite{BruchhaeuserBR12} and set it to 1.6072872.
The flow is driven by the Dirichlet inflow profile \(\boldsymbol v=\boldsymbol v_D\) on the left boundary
$\Gamma_{\mathrm{in}} \coloneq \{0\}\times[0,H],$ given by
\begin{equation}
\label{eq:inflow}
\boldsymbol v_D(\boldsymbol{x},t)
=
\begin{pmatrix}
\sin\!\left(\frac{\pi}{8}t\right)\dfrac{6y(H-y)}{H^2}\\[0.3em]
0
\end{pmatrix}
\qquad \text{on }\Gamma_{\mathrm{in}},
\end{equation}
where \(H=0.41\) denotes the channel height. On the wall boundary \(\Gamma_{\mathrm{wall}}\) and on the obstacle boundary we prescribe no-slip boundary conditions \(\boldsymbol v=\boldsymbol 0\). On the outflow boundary \(\Gamma_{\mathrm{out}}\) we impose a do-nothing condition~\cite{BruchhaeuserHRT96}.

In Table~\ref{tab:dfg-drag-all} we present the results of an adaptive refinement run together with a comparison of the algebraic performance of a GMG-preconditioned Newton-Krylov method and an ILU method using a dG($1$)-$Q_2/P_1^{\mathrm{disc}}$ discretization.
The adaptive refinement reduces the mean drag error \eqref{eq:Jmeandrag} monotonically. The drag value converges to the reference value and the effectivity indices stay close to one on all loops. This shows that the DWR estimator provides a reliable approximation of the goal error.
The observed decay of the goal error with respect to the space-time number of degrees of freedom is at least linear and becomes higher on the finer loops. Measured against total algebraic work, the same trend is obtained more efficiently with GMG than with ILU.
\begin{table}[H]
\centering
\caption{Adaptive refinement and algebraic performance for the
dG($1$)-$Q_2/P_1^{\mathrm{disc}}$ discretization with goal functional
$J_{\mathrm{drag}}$. \subref{tab:dfg-drag-a} reports the mean-drag value, the
goal error, the space-time estimator contributions, and the effectivity index.
\subref{tab:dfg-drag-b} compares the GMG and ILU preconditioners in terms of
nonlinear iterations per slab, primal linear iterations per nonlinear step,
total algebraic work, and the resulting convergence rates with respect to work.}
\label{tab:dfg-drag-all}
\begin{subcaptionblock}{\textwidth}
\caption{Goal quantity, error estimator, and effectivity.}
\label{tab:dfg-drag-a}
\centering
\resizebox{.99\linewidth}{!}{%
\begin{tabular}{
S[scientific-notation=false,round-precision=0,table-format=1]
S[scientific-notation=false,round-precision=0,table-format=3]
S[scientific-notation=false,round-precision=0,table-format=7]
S[table-format=1.3,round-precision=3]
S[table-format=1.2e-1]
S[table-format=1.2e-1]
S[table-format=1.2e-1]
S[table-format=1.2e-1]
S[table-format=1.4]}
\toprule
\mc{$\ell$} & \mc{$N_{\text{time}}$}
& \mc{$N_{\text{tot}}$}
& \mc{$J_{\text{drag}}(u_{\tau h})$}
& \mc{$\|e_{\text{drag}}\|$}
& \mc{$\tilde{\eta}_{h}$} & \mc{$\tilde{\eta}_\tau$}
& \mc{$\tilde{\eta}_{\tau h}$} & \mc{$\Ieff$}
\\
\midrule
1 & 160 &  606720 & 1.389 & 2.137e-01 & 1.565e-01 & 9.082e-02 & 2.473e-01 & 1.1571 \\
2 & 176 & 1341824 & 1.516 & 8.740e-02 & 5.680e-02 & 3.951e-02 & 9.631e-02 & 1.1020 \\
3 & 193 & 3004238 & 1.584 & 1.930e-02 & 1.346e-02 & 8.340e-03 & 2.180e-02 & 1.1294 \\
4 & 212 & 7023508 & 1.601 & 1.860e-03 & 1.222e-03 & 8.710e-04 & 2.093e-03 & 1.1250 \\
\bottomrule
\end{tabular}
}
\end{subcaptionblock}

\medskip

\begin{subcaptionblock}{\textwidth}
\caption{Solver robustness, total work, and convergence rates.}
\label{tab:dfg-drag-b}
\centering
\resizebox{.99\linewidth}{!}{%
\begin{tabular}{
S[scientific-notation=false,round-precision=0,table-format=1]
S[table-format=2.2]
S[table-format=2.2]
S[table-format=2.2]
S[table-format=2.2]
S[table-format=1.2e2]
S[table-format=1.2e2]
S[table-format=1.2]
S[table-format=1.2]
S[table-format=1.2]}
\toprule
&
\multicolumn{2}{c}{$n_{\mathrm{NL}}/N_{\mathrm{time}}$} &
\multicolumn{2}{c}{$n_{\mathrm{L}}^{p}/n_{\mathrm{NL}}$} &
\multicolumn{2}{c}{$W_{\mathrm{tot}}$} &
\multicolumn{3}{c}{rates of $\|e_{\mathrm{drag}}\|$}
\\
\cmidrule(lr){2-3}\cmidrule(lr){4-5}\cmidrule(lr){6-7}\cmidrule(lr){8-10}
\mc{$\ell$}
& \mc{GMG} & \mc{ILU}
& \mc{GMG} & \mc{ILU}
& \mc{GMG} & \mc{ILU}
& \mc{$N_{\mathrm{tot}}$}
& \mc{$W^{\mathrm{GMG}}$}
& \mc{$W^{\mathrm{ILU}}$}
\\
\midrule
1 & 4.3937 &  8.8719 &   1.0000 & 111.5196 & 3.5076e+07 & 6.3608e+08 & \multicolumn{1}{c}{--} & \multicolumn{1}{c}{--} & \multicolumn{1}{c}{--} \\
2 & 4.2045 &  9.2500 &  10.2622 & 112.0510 & 1.3350e+08 & 1.4861e+09 & 1.126 & 0.669 & 1.054 \\
3 & 4.2539 &  6.4197 &  12.1133 & 109.3930 & 3.2757e+08 & 2.4019e+09 & 1.874 & 1.683 & 3.146 \\
4 & 4.6368 &  6.8538 &  13.2970 & 105.3066 & 8.3955e+08 & 3.7967e+09 & 2.755 & 2.486 & 5.110 \\
\bottomrule
\end{tabular}
}
\end{subcaptionblock}
\end{table}
This is consistent with the solver data. The GMG-preconditioned Newton-Krylov method remains robust over the full adaptive sequence, whereas ILU requires
substantially more nonlinear iterations and much larger primal linear iteration
counts. The larger number of nonlinear iterations for ILU can be explained by
the lower accuracy of the linear solves in the Newton steps. Since the resulting
Newton corrections are less accurate than with GMG, the convergence of the outer
Newton iteration is correspondingly weaker. Overall, the results indicate that
the adaptive strategy is effective for the drag goal quantity \eqref{eq:Jmeandrag} and that GMG is the more robust and efficient preconditioner for this problem.

\begin{figure}[H]
  \centering
\resizebox{.99\linewidth}{!}{%
\includegraphics{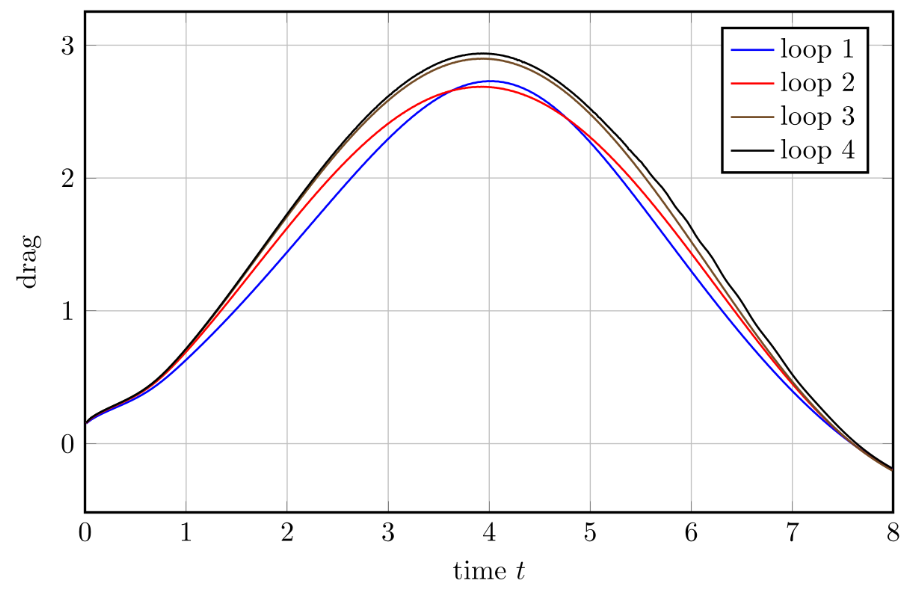}
  }
  \caption{Drag over all DWR loops}
  \label{fig:drag-lift-all}
\end{figure}

\begin{figure}[H]
\centering

\begin{subcaptionblock}{0.95\textwidth}
  \centering
  \includegraphics[width=\textwidth]{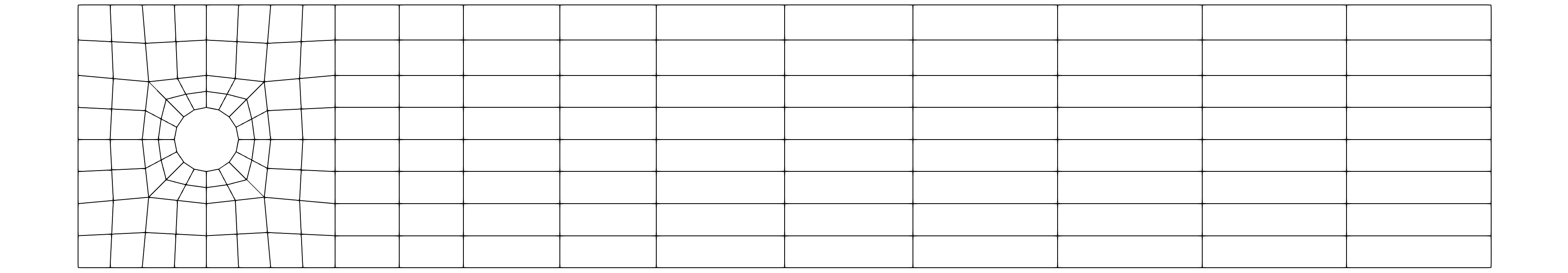}
  \caption{Adaptive mesh after loop 1.}
\end{subcaptionblock}
\hfill
\begin{subcaptionblock}{0.95\textwidth}
  \centering
  \includegraphics[width=\textwidth]{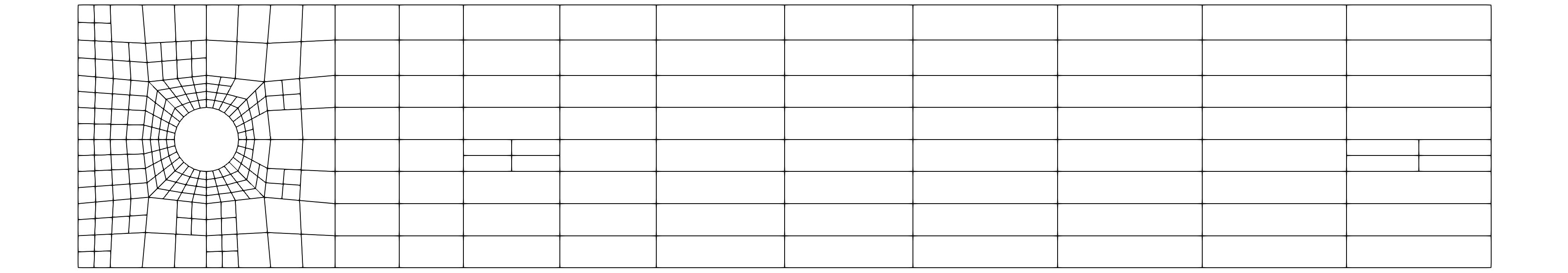}
  \caption{Adaptive mesh after loop 2.}
\end{subcaptionblock}
\hfill
\begin{subcaptionblock}{0.95\textwidth}
  \centering
  \includegraphics[width=\textwidth]{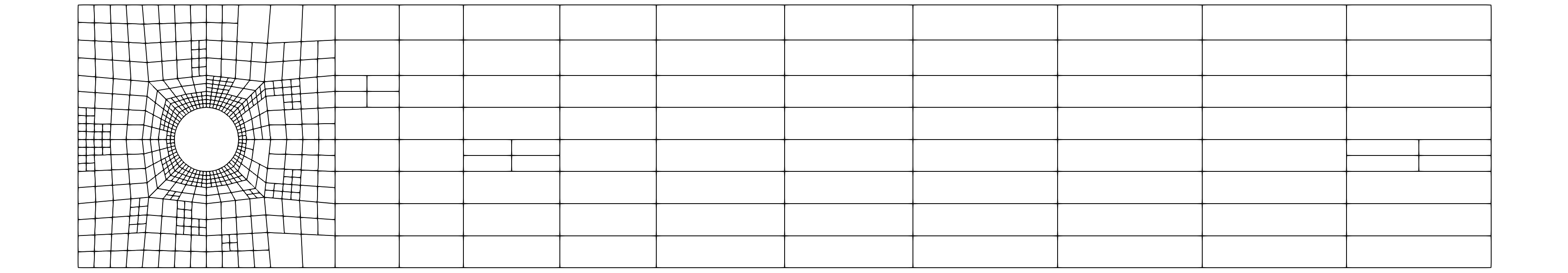}
  \caption{Adaptive mesh after loop 3.}
\end{subcaptionblock}
\begin{subcaptionblock}{0.95\textwidth}
  \centering
  \includegraphics[width=\textwidth]{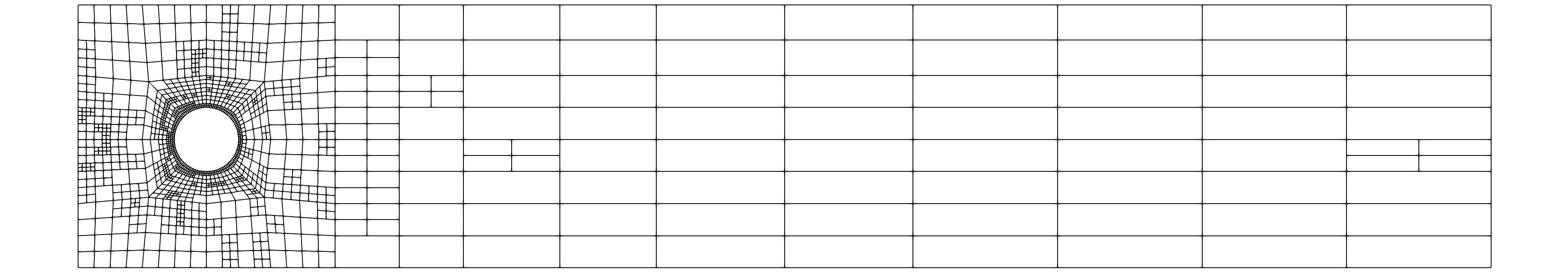}
  \caption{Adaptive mesh after loop 4.}
\end{subcaptionblock}
\hfill
\begin{subcaptionblock}{0.95\textwidth}
  \centering
  \includegraphics[width=\textwidth]{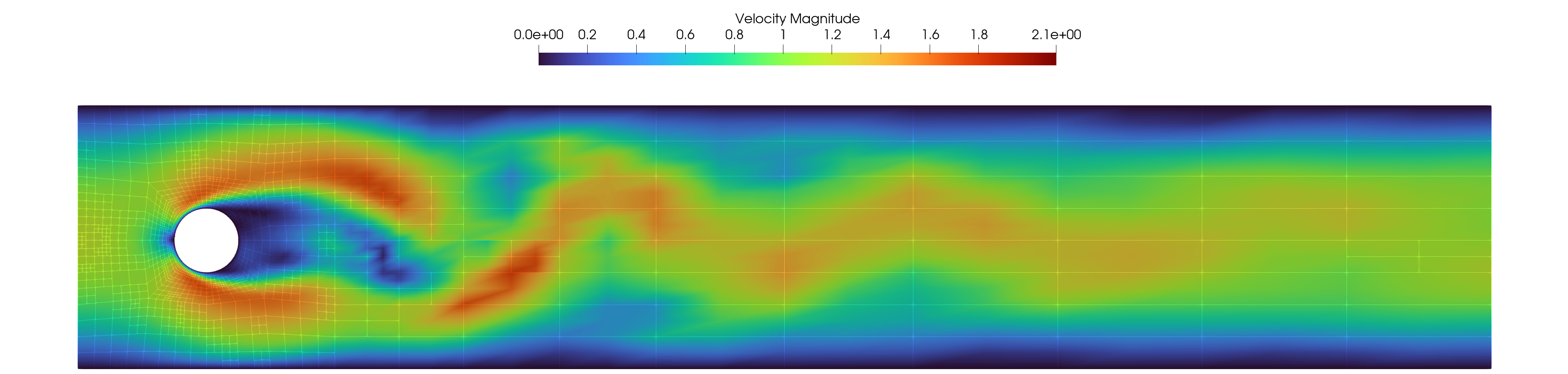}
  \caption{Discrete solution on the final loop.}
\end{subcaptionblock}

\caption{Adaptive refinement history and final-loop solution at $t=5$. The first four panels show the meshes generated over the adaptive loops. The last two panels show a refined view of the final mesh and the corresponding discrete solution.}
\label{fig:adaptive-meshes-solution}
\end{figure}

The temporal evolution of the drag coefficient is compared for the different DWR refinement loops in Fig.~\ref{fig:drag-lift-all}. We note that the curves converge by increasing the adaptive refinement (cf. loop 3 and 4) and are in good comparison to the results given in \cite[Fig.~10]{BruchhaeuserBR12} and \cite[Fig.~6(a)]{RoThiKoeWi23}, respectively.
Finally, in Fig.~\ref{fig:adaptive-meshes-solution} we present the discrete primal solution on the final DWR loop as well as the adaptively refined spatial meshes for all four DWR loops at time point $t=5$.
The spatial refinement remains concentrated close to the cylinder and in the wake, as confirmed by the selected meshes and solution snapshots. The drag and lift histories are smooth and consistent with the benchmark results, which indicates that the adaptive discretization resolves the relevant time dependent features and the Stokes projections avoids spurious oscillations.

\section{Conclusion}
\label{sec:6:conclusion}
We have developed a goal-oriented space-time adaptive finite element method for the incompressible Navier--Stokes equations based on the Dual Weighted Residual method. The approach combines discontinuous Galerkin time discretization, inf-sup stable mixed finite elements in space, adaptive refinement in space and time, and a slab-wise Newton--Krylov solver with geometric multigrid preconditioning. The numerical results show that this combination yields effective control of the goal quantities on dynamically adapted meshes.

The main computational benefit of the DWR approach is its problem-driven refinement. The estimator identifies those regions in space and time that are most relevant for the chosen output and thereby focuses resolution where it matters for the quantity of interest. At the same time, the results show that this accuracy comes at a substantial cost. In addition to the primal problem, the adaptive loop requires the repeated solution of a dual problem and the evaluation of weighted residual contributions. The efficiency of the overall approach therefore depends not only on the quality of the estimator, but equally on the robustness of the underlying iterative solvers.

The Newton iteration counts are essentially unaffected by the choice of linear preconditioner, whereas the GMRES iteration counts depend strongly on it. Thus, the additional cost of the DWR framework is governed primarily by the linear solver performance for the primal and dual subproblems. In this setting, geometric multigrid provides the required robustness, while ILU does not yield a competitive performance. This is especially important for higher-order discretizations, where the linear systems become more demanding and the dual solve adds a further substantial burden. For goal-oriented space-time adaptivity, robust preconditioning is therefore essential.

Overall, the results show that DWR-based space-time adaptivity is effective for goal-oriented error control, provided that it is paired with scalable linear algebra. The practical cost of the method is determined by the interplay between estimator quality and solver performance. Future work should therefore address both aspects together, in particular through cheaper dual solves, more selective dual updates, and improved multilevel solvers for higher-order and larger-scale space-time discretizations.

\FloatBarrier
\section*{Acknowledgement}
\noindent
Computational resources (HPC cluster HSUper) have been provided by the project hpc.bw, funded by dtec.bw - Digitalization and Technology Research Center of the Bundeswehr. dtec.bw is funded by the European Union - NextGenerationEU.

\bibliographystyle{plain}
\bibliography{bruchhaeuser-refs}

@article{RoThiKoeWi23,
  title = {Tensor-{P}roduct {S}pace-{T}ime {G}oal-{O}riented {E}rror {C}ontrol and {A}daptivity {W}ith {P}artition-of-{U}nity {D}ual-{W}eighted {R}esiduals for {N}onstationary {F}low {P}roblems},
  author = {Julian Roth and Jan Philipp Thiele and Uwe Köcher and Thomas Wick},
  pages = {185--214},
  volume = {24},
  number = {1},
  journal = {Comput. Methods Appl. Math.},
  doi = {10.1515/cmam-2022-0200},
  year = {2024},
  lastchecked = {2024}
}

@article{besier_pressure_2012,
	title = {On the pressure approximation in nonstationary incompressible flow simulations on dynamically varying spatial meshes: {PRESSURE} {APPROXIMATION} {ON} {DYNAMIC} {MESHES}},
	volume = {69},
	issn = {02712091},
	url = {https://onlinelibrary.wiley.com/doi/10.1002/fld.2625},
	doi = {10.1002/fld.2625},
	shorttitle = {On the pressure approximation in nonstationary incompressible flow simulations on dynamically varying spatial meshes},
	pages = {1045--1064},
	number = {6},
	journaltitle = {International Journal for Numerical Methods in Fluids},
	shortjournal = {Int. J. Numer. Meth. Fluids},
	author = {Besier, Michael and Wollner, Winnifried},
	year = {2012},
}

@ARTICLE{BruchhaeuserABBFGHHKKMMPPSTWZ24,
  title = {The deal.{II} library, Version 9.6},
  volume = {32},
  issn = {1569-3953},
  doi = {10.1515/jnma-2024-0137},
  pages = {369--380},
  number = {4}, 
  journal = {Journal of Numerical Mathematics},
  author = {Africa, Pasquale C. and Arndt, Daniel and Bangerth, Wolfgang and Blais, Bruno and Fehling, Marc and Gassmöller, Rene and Heister, Timo and Heltai, Luca and Kinnewig, Sebastian and Kronbichler, Martin and Maier, Matthias and Munch, Peter and Schreter-Fleischhacker, Magdalena and Thiele, Jan P. and Turcksin, Bruno and Wells, David and Yushutin, Vladimir},
  year ={2024}
}

@ARTICLE{BruchhaeuserAB23,
  title={A geometric multigrid method for space-time finite element discretizations of the Navier--Stokes equations and its application to 3d flow simulation},
  author={Anselmann, Mathias and Bause, Markus},
  journal={ACM Transactions on Mathematical Software},
  volume={49},
  number={1},
  pages={1--25},
  year={2023},
  publisher={ACM New York, NY}
}

@ARTICLE{BruchhaeuserBBEMTW25,
      title={Anisotropic space-time goal-oriented error control and mesh adaptivity for convection-diffusion-reaction equations},
      author={Bause, M. and Bruchhäuser, M. P. and Endtmayer, B. and Margenberg, N. and Toulopoulos, I. and Wick, T.},
      journal = {J. Numer. Math},
      year={2025},
      archivePrefix={arXiv},
      primaryClass={math.NA},
      doi={10.48550/arXiv.2504.04951}
}

@BOOK{BruchhaeuserBR03,
 author={W. Bangerth and R. Rannacher},
 title={Adaptive Finite Element Methods for Differential Equations},
 address={Basel},
 OPTseries = {Lectures in Mathematics, ETH Z\"urich},
publisher =    {Birkh{\"a}user Verlag},
 year={2003}
}

@ARTICLE{BruchhaeuserBGR10,
url = {https://doi.org/10.2478/cmam-2010-0001},
title = {Adaptive Galerkin Finite Element Methods for the Wave Equation},
author = {W. Bangerth and M. Geiger and R. Rannacher},
pages = {3--48},
volume = {10},
number = {1},
journal = {Computational Methods in Applied Mathematics},
doi = {doi:10.2478/cmam-2010-0001},
year = {2010}
}

@ARTICLE{BruchhaeuserBR01,
    AUTHOR = {Becker, R. and Rannacher, R.},
     TITLE = {An optimal control approach to a posteriori error estimation
              in finite element methods},
   JOURNAL = {Acta Numer.},
  FJOURNAL = {Acta Numerica},
    VOLUME = {10},
      YEAR = {2001},
     PAGES = {1--102},
      ISSN = {0962-4929},
	doi = {10.1017/S0962492901000010}
}

@ARTICLE{BruchhaeuserBR12,
  author = {Besier, M. and Rannacher, R.},
  title = {Goal-oriented space-time adaptivity in the finite element {G}alerkin method for the computation of nonstationary incompressible flow},
  journal = {Int. J. Num. Methods Fluids},
  volume = {70},
  number = {9},
  pages = {1139--1166},
  year = {2012},
  doi = {10.1002/fld.2735}
}

@article{10.1002/fld.2625,
author = {Besier, Michael and Wollner, Winnifried},
title = {On the pressure approximation in nonstationary incompressible flow simulations on dynamically varying spatial meshes},
journal = {International Journal for Numerical Methods in Fluids},
volume = {69},
number = {6},
pages = {1045-1064},
doi = {https://doi.org/10.1002/fld.2625},
url = {https://onlinelibrary.wiley.com/doi/abs/10.1002/fld.2625},
eprint = {https://onlinelibrary.wiley.com/doi/pdf/10.1002/fld.2625},
year = {2012}
}

@INCOLLECTION{RePEc:spr:978-3-642-18775,
title = {Adaptive Computation of Reactive Flows with Local Mesh Refinement and Model Adaptation},
author = {Braack, Malte and Ern, Alexandre},
year = {2004},
pages = {159-168},
publisher = {Springer},
url = {https://EconPapers.repec.org/RePEc:spr:sprchp:978-3-642-18775-9_13}
}

@article{10.1093/imanum/drt001,
    author = {Brenner, Andreas and Bänsch, Eberhard and Bause, Markus},
    title = {A priori error analysis for finite element approximations of the Stokes problem on dynamic meshes},
    journal = {IMA Journal of Numerical Analysis},
    volume = {34},
    number = {1},
    pages = {123-146},
    year = {2013},
    month = {05},
    issn = {0272-4979},
    doi = {10.1093/imanum/drt001},
    url = {https://doi.org/10.1093/imanum/drt001},
    eprint = {https://academic.oup.com/imajna/article-pdf/34/1/123/1950435/drt001.pdf},
}

@InProceedings{BruchhaeuserBB24,
author={Bruchhäuser, Marius Paul and Bause, Markus},
editor={Sequeira, Ad{\'e}lia and Silvestre, Ana and Valtchev, Svilen S. and Janela, Jo{\~a}o},
title={Numerical Study of Approximation Techniques for the Temporal Weights to the DWR Method},
booktitle={Numerical Mathematics and Advanced Applications ENUMATH 2023, Volume 1},
date={2025},
publisher={Springer Nature Switzerland},
address={Cham},
pages={177--187},
doi = {10.1007/978-3-031-86173-4\_18}
}

@ARTICLE{BruchhaeuserBKB22,
title = {On the Implementation of an Adaptive Multirate Framework for Coupled Transport and Flow},
journal = {Journal of Scientific Computing},
volume = {93},
number = {59},
date = {2022},
doi = {10.1007/s10915-022-02026-z},
author = {Bruchhäuser, Marius Paul and Köcher, Uwe and Bause, Markus}
}

@phdthesis{BruchhaeuserB22,
  author = {Bruchh\"auser, M. P.},
  title = {Goal-oriented space-time adaptivity for a multirate approach to coupled flow and transport},
  school = {Helmut-Schmidt-University/University of the German Federal Armed Forces Hamburg},
  year = {2022},
  doi = {10.24405/14380}
}

@BOOK{BruchhaeuserCO84,
	author    = {Carey, G.F. and Oden, J.T.},
	title     = {Finite elements, computational aspects, Vol. III},
	series    = {The Texas finite element series},
	publisher = {Prentice-Hall},
	address   = {Englewood Cliffs, New Jersey},
	year      = {1984}
}

@article{ern-2017,
author = {Ern, Alexandre and Smears, Iain and Vohral\'{\i}k, Martin},
title = {Guaranteed, Locally Space-Time Efficient, and Polynomial-Degree Robust a Posteriori Error Estimates for High-Order Discretizations of Parabolic Problems},
journal = {SIAM Journal on Numerical Analysis},
volume = {55},
number = {6},
pages = {2811-2834},
year = {2017},
doi = {10.1137/16M1097626}
}

@article{falgout-2017,
	author = {Falgout, R. D. and Friedhoff, S. and Kolev, Tz. V. and MacLachlan, S. P. and Schroder, J. B. and Vandewalle, S.},
	journal = {Computing and Visualization in Science},
	month = {8},
	number = {4-5},
	pages = {123--143},
	title = {{Multigrid methods with space–time concurrency}},
	volume = {18},
	year = {2017},
	doi = {10.1007/s00791-017-0283-9},
	url = {https://doi.org/10.1007/s00791-017-0283-9},
}

@article{10.1137/18M1172466,
author = {Neum\"{u}ller, Martin and Smears, Iain},
title = {Time-Parallel Iterative Solvers for Parabolic Evolution Equations},
journal = {SIAM Journal on Scientific Computing},
volume = {41},
number = {1},
pages = {C28-C51},
year = {2019},
doi = {10.1137/18M1172466},
URL = {  https://doi.org/10.1137/18M1172466},
eprint = { 
        https://doi.org/10.1137/18M1172466
}
}

@article{gmeiner-2015,
	author = {Gmeiner, Björn and Rüde, Ulrich and Stengel, Holger and Waluga, Christian and Wohlmuth, Barbara},
	journal = {Numerical Mathematics Theory Methods and Applications},
	month = {2},
	number = {1},
	pages = {22--46},
	title = {{Towards Textbook Efficiency for Parallel Multigrid}},
	volume = {8},
	year = {2015},
	doi = {10.4208/nmtma.2015.w10si},
	url = {https://doi.org/10.4208/nmtma.2015.w10si},
}

@book{hackbusch1985multigrid,
  title={Multi-Grid Methods and Applications},
  author={Hackbusch, Wolfgang},
  year={1985},
  publisher={Springer-Verlag},
  address={Berlin Heidelberg New York Tokyo},
  series={Springer Series in Computational Mathematics},
  volume={4},
  isbn={3-540-12761-5},
  doi={10.1007/978-3-662-02427-0}
}

@article{BruchhaeuserHR82,
author = {Heywood, John G. and Rannacher, Rolf},
title = {Finite Element Approximation of the Nonstationary Navier–Stokes Problem. I. Regularity of Solutions and Second-Order Error Estimates for Spatial Discretization},
journal = {SIAM Journal on Numerical Analysis},
volume = {19},
number = {2},
pages = {275-311},
year = {1982},
doi = {10.1137/0719018},
URL = {https://doi.org/10.1137/071901},
eprint = {https://doi.org/10.1137/0719018}
}

@ARTICLE{BruchhaeuserHRT96,
  title={Artificial boundaries and flux and pressure conditions for the incompressible Navier--Stokes equations},
  author={Heywood, John G and Rannacher, Rolf and Turek, Stefan},
  journal={International Journal for numerical methods in fluids},
  volume={22},
  number={5},
  pages={325--352},
  year={1996},
  publisher={Wiley Online Library}
}

@article{hussain-2012,
	author = {Hussain, S.},
	journal = {The Open Numerical Methods Journal},
	month = {1},
	number = {1},
	pages = {35--45},
	title = {{A Note on Accurate and Efficient Higher Order Galerkin Time Stepping Schemes for the Nonstationary Stokes Equations}},
	volume = {4},
	year = {2012},
	doi = {10.2174/1876389801204010035},
	url = {https://doi.org/10.2174/1876389801204010035},
}

@BOOK{BruchhaeuserJ16,
  author = 	{John, Volker},
  title = 	{Finite Element Methods for Incompressible Flow Problems},
  series = 	{Springer series in computational mathematics},
  year = 	{2016},
  edition = 	{1. Edition},
  publisher = 	{Springer Cham},
  address = 	{Berlin; Heidelberg},
  volume = 	{51},
  isbn = 	{9783319457505},
  issn = 	{2198-3712},
  doi = 	{10.1007/978-3-319-45750-5},
  language = 	{en}
}

@article{10.1002/fld.195,
author = {John, Volker and Matthies, Gunar},
title = {Higher-order finite element discretizations in a benchmark problem for incompressible flows},
journal = {International Journal for Numerical Methods in Fluids},
volume = {37},
number = {8},
pages = {885-903},
doi = {https://doi.org/10.1002/fld.195},
url = {https://onlinelibrary.wiley.com/doi/abs/10.1002/fld.195},
eprint = {https://onlinelibrary.wiley.com/doi/pdf/10.1002/fld.195},
year = {2001}
}

@article{langer-2021,
	author = {Langer, Ulrich and Schafelner, Andreas},
	journal = {Journal of Numerical Mathematics},
	month = {11},
	number = {4},
	pages = {247--266},
	title = {{Adaptive space–time finite element methods for parabolic optimal control problems}},
	volume = {30},
	year = {2021},
	doi = {10.1515/jnma-2021-0059},
	url = {https://doi.org/10.1515/jnma-2021-0059},
}

@book{10.1515/9783110548488,
url = {https://doi.org/10.1515/9783110548488},
title = {Space-Time Methods},
title = {Applications to Partial Differential Equations},
editor = {Ulrich Langer and Olaf Steinbach},
publisher = {De Gruyter},
address = {Berlin, Boston},
doi = {doi:10.1515/9783110548488},
isbn = {9783110548488},
year = {2019},
lastchecked = {2026-03-25}
}

@ARTICLE{BruchhaeuserMT02,
  title={The Inf-Sup Condition for the Mapped Q k- P k-1 disc Element in Arbitrary Space Dimensions},
  author={Matthies, G. and Tobiska, L.},
  journal={Computing},
  volume={69},
  number={2},
  pages={119--139},
  year={2002},
  publisher={Springer}
}

@inproceedings{BruchhaeuserN71,
  title={{\"U}ber ein Variationsprinzip zur L{\"o}sung von Dirichlet-Problemen bei Verwendung von Teilr{\"a}umen, die keinen Randbedingungen unterworfen sind},
  author={Nitsche, Joachim},
  booktitle={Abhandlungen aus dem mathematischen Seminar der Universit{\"a}t Hamburg},
  volume={36},
  number={1},
  pages={9--15},
  year={1971},
  organization={Springer}
}

@BOOK{BruchhaeuserR17,
  title={Fluid-structure interactions: models, analysis and finite elements},
  author={Richter, Thomas},
  volume={118},
  year={2017},
  publisher={Springer}
}

@ARTICLE{BruchhaeuserSV08,
  author = {Schmich, M. and Vexler, B.},
  title = {Adaptivity with dynamic meshes for space-time finite element discretizations of parabolic equations},
  journal = {SIAM J. Sci. Comput.},
  volume = {30},
  number = {1},
  pages = {369--393},
  year = {2008},
  doi = {10.1137/060670468}
}

@incollection{BruchhaeuserST96,
  title={Benchmark computations of laminar flow around a cylinder},
  author={Sch{\"a}fer, Michael and Turek, Stefan and Durst, Franz and Krause, Egon and Rannacher, Rolf},
  booktitle={Flow simulation with high-performance computers II: DFG priority research programme results 1993--1995},
  pages={547--566},
  year={1996},
  publisher={Springer}
}

@BOOK{KelleyIterativeMethodsLinear1995,
  title={Iterative methods for linear and nonlinear equations},
  author={Kelley, Carl T},
  year={1995},
  publisher={SIAM}
}

@article{25M1734142,
author = {Margenberg, Nils and Bause, Markus and Munch, Peter},
title = {An \(hp\) Multigrid Approach for Tensor-Product Space-Time Finite Element Discretizations of the Stokes Equations},
journal = {SIAM Journal on Scientific Computing},
volume = {47},
number = {6},
pages = {B1503-B1529},
year = {2025},
doi = {10.1137/25M1734142},
URL = { 
        https://doi.org/10.1137/25M1734142
},
eprint = { 
        https://doi.org/10.1137/25M1734142
}
}

@book{10.1137/1.9781611971057,
author = {McCormick, Stephen F.},
editor = {Stephen F. McCormick},
title = {Multigrid Methods},
publisher = {Society for Industrial and Applied Mathematics},
year = {1987},
doi = {10.1137/1.9781611971057},
address = {},
edition   = {},
}

@BOOK{Saad1996,
  title={Iterative methods for sparse linear systems},
  author={Saad, Yousef},
  year={2003},
  publisher={SIAM}
}

@InCollection{SteinbachYang2019,
  author    = {Steinbach, Olaf and Yang, Huidong},
  title     = {Space-time finite element methods for parabolic evolution equations: discretization, a posteriori error estimation, adaptivity and solution},
  booktitle = {Space-time methods: Applications to partial differential equations},
  year      = {2019},
  editor    = {Taghavy, H. and others},
  series    = {De Gruyter},
  pages     = {207--248},
  address   = {Berlin},
  publisher = {De Gruyter},
  doi       = {10.1515/9783110548488-007}
}

@book{steinbach-2018,
	author = {Steinbach, Olaf and Yang, Huidong},
	booktitle = {Lecture notes in computer science},
	month = {1},
	pages = {66--73},
	title = {{An Algebraic Multigrid Method for an Adaptive Space–Time Finite Element Discretization}},
	year = {2018},
	doi = {10.1007/978-3-319-73441-5\_6},
}

@book{1977iii,
title = {Navier--Stokes Equations},
author = {Roger Temam},
publisher = {North-Holland},
year = {1984},
isbn = {0 444 87558},
address = {Amsterdam},
}

@book{10.5555/1211262, 
author = {Thom\'{e}e, Vidar}, 
title = {Galerkin Finite Element Methods for Parabolic Problems (Springer Series in Computational Mathematics)}, 
year = {2006}, 
isbn = {3540331212}, 
publisher = {Springer-Verlag}, 
address = {Berlin, Heidelberg} 
}

@ARTICLE{Vanka1985,
  title={Block-implicit multigrid solution of Navier-Stokes equations in primitive variables},
  author={Vanka, S Pratap},
  journal={Journal of Computational Physics},
  volume={65},
  number={1},
  pages={138--158},
  year={1986},
  publisher={Elsevier},
  doi = {10.1016/0021-9991(86)90008-2}
}

@article{WIENERS2023294,
title = {A space-time discontinuous Galerkin discretization for the linear transport equation},
journal = {Computers \& Mathematics with Applications},
volume = {152},
pages = {294-307},
year = {2023},
issn = {0898-1221},
doi = {https://doi.org/10.1016/j.camwa.2023.10.031},
url = {https://www.sciencedirect.com/science/article/pii/S0898122123004789},
author = {Christian Wieners}
}

@article{BruchhaeuserWW21,
  title={Optimization with nonstationary, nonlinear monolithic fluid-structure interaction},
  author={Wick, Thomas and Wollner, Winnifried},
  journal={International Journal for Numerical Methods in Engineering},
  volume={122},
  number={19},
  pages={5430--5449},
  year={2021},
  publisher={Wiley Online Library}
}

\end{document}